\documentclass[dvips,preprint,aap]{imsart}

\usepackage{imsart}
\usepackage{pstricks} 
\usepackage{graphicx} 
\usepackage{amsmath}
\usepackage{amsfonts}
\usepackage{amssymb}
\usepackage{enumerate}
\usepackage{ulem}
\usepackage{version}
\usepackage{verbatim}
\usepackage{moreverb}
\usepackage{float}

\usepackage{amsthm}

\usepackage[english]{babel} 
\usepackage[latin1]{inputenc}
\usepackage[T1]{fontenc}
\usepackage{pst-all}
\usepackage{dsfont}
\RequirePackage[OT1]{fontenc}
\RequirePackage{amsmath}
\RequirePackage[numbers]{natbib}
\RequirePackage[colorlinks,citecolor=blue,urlcolor=blue]{hyperref}

\usepackage{yfonts}

\startlocaldefs
\newtheoremstyle{them}
{5mm}
{3pt}
{\itshape}
{}
{\bfseries}
{ :}
{\newline}
{}

\newtheoremstyle{exdef}
{5mm}
{3pt}
{}
{}
{\bfseries}
{:}
{\newline}
{}
\theoremstyle{them}
\newtheorem{thm}{Theorem}[section]
\newtheorem{prop}[thm]{Proposition}

\newtheorem{cor}[thm]{Corollary}

\newtheorem{lem}[thm]{Lemma}

\newtheorem{hyp}[thm]{Hypothesis}
\newtheorem{sthyp}[thm]{Standing Hypothesis}

\theoremstyle{exdef}
\newtheorem{defin}[thm]{Definition}
\newtheorem{ex}{Example}

\newcommand{\R}{\mathbb{R}}
\newcommand{\N}{\mathbb{N}}
\newcommand{\C}{\mathbb{C}}

\renewcommand{\L}{\mathcal{L}}
\newcommand{\Z}{\mathbb{Z}}

\newcommand{\E}{\mathbb{E}}
\newcommand{\p}{\mathbb{P}}
\newcommand{\e}{e^}
\renewcommand{\bar}{\overline}

\newcommand{\ep}{\varepsilon}
\newcommand{\iD}{\mathring{\Delta}}
\newcommand{\adh}{\overline}

\newcommand{\ap}{\leadsto_{ap}}
\def\1{{\rm 1\mskip-4.4mu l}} 
\def\ds{\displaystyle}
\def\restriction#1#2{\mathchoice
              {\setbox1\hbox{${\displaystyle #1}_{\scriptstyle #2}$}
              \restrictionaux{#1}{#2}}
              {\setbox1\hbox{${\textstyle #1}_{\scriptstyle #2}$}
              \restrictionaux{#1}{#2}}
              {\setbox1\hbox{${\scriptstyle #1}_{\scriptscriptstyle #2}$}
              \restrictionaux{#1}{#2}}
              {\setbox1\hbox{${\scriptscriptstyle #1}_{\scriptscriptstyle #2}$}
              \restrictionaux{#1}{#2}}}
\def\restrictionaux#1#2{{#1\,\smash{\vrule height .8\ht1 depth .85\dp1}}_{\,#2}} 
\newenvironment{pre}{ \noindent \textbf{\underline{Proof :\newline}} }{ \begin{flushright} $\blacksquare$ \end{flushright}}

\endlocaldefs

\begin{document}

\begin{frontmatter}

\title{Quasi-stationary distributions for stochastic approximation algorithms with constant step size}
\runtitle{QSD for stochastic approximation algorithms}


\begin{aug}
\author{\fnms{Bastien} \snm{Marmet}\ead[label=e1]{bastien.marmet@unine.ch}\thanksref{t1}}
\affiliation{University of Neuchatel}
\address{Institut de Math\'ematiques, Universit\'e de Neuch\^atel,\\ Rue Emile-Argand 11. Neuch\^atel. Switzerland.\\ \printead{e1}}

\thankstext{t1}{The author acknowledge financial support from the Swiss National Foundation Grant FN
200021-138242/1}

\runauthor{Bastien Marmet}
\end{aug}
\begin{abstract}
: In this paper we investigate quasi-stationary distributions $\mu^N$ of stochastic approximation algorithms with constant step size which can be viewed as random perturbations of a time-continuous dynamical system. Inspired by ecological models these processes have a closed absorbing set corresponding to extinction. Under some large deviation assumptions and the existence of an interior attractor for the ODE, we show that the weak* limit points of the QSD $\mu^N$ are invariant measures for the ODE with support in the interior attractors. 
\end{abstract}



\end{frontmatter}


\tableofcontents

\section{Introduction}
\hspace{1mm}\newline
One of the most considered issue in theoretical ecology is to find out under what kind of conditions one can expect a population of interacting species (animals, plants, microorganisms, ...) to survive on the long term with no extinctions. When these conditions are met the interacting populations are said to persist or coexist. In the past, differential equations and nonlinear difference equations have been used to model these phenomena. Famous examples are Lotka \cite{lotka-25} and Volterra \cite{volterra-26} work on competitive and predator-prey interactions, Thompson \cite{thompson-24}, Nicholson and Bailey \cite{nicholson-bailey-35} on host-parasite interactions, and Kermack and McKendrick \cite{kermack-mckendrick-27} on disease outbreaks.
For these deterministic models, persistence definitions sometimes vary but most authors link persistence with the existence of an attractor bounded away from the extinction states, in which case persistence holds over an infinite time horizon, see e.g. \cite{Sch06}. In order to refine these models and allow for some "roughness" and/or influence of unpredictable outer events, randomness has been added to these models, leading to Markov processes models. However, extinctions being absorbent states and species dying out with positive probabilities, the underlying theory of Markov processes shows that, in finite time, extinction is inevitable. Yet, in the real world, with large sized pools of population, we don't observe that inevitable extinction. This finite extinction time may then be very large and the system may remain in some sort of "metastable state" bounded away from extinction for a long time.
In \cite{FauSch11}, Faure and Schreiber studied this problem for randomly perturbed discrete time dynamical systems, showing that, under the appropriate assumptions about the random perturbations, there exists a positive attractor (i.e. an attractor which is bounded away from extinction states) for the unperturbed system, which implies two things as the number of individuals or particles gets large. First, when they exist, quasi-stationary distributions concentrate on the positive attractors of the unperturbed system. Second, the expected time to extinction for systems starting according to this quasi-stationary distribution grows exponentially with the system size. The aim of this paper is to extend their approach to a class of discrete time Markov process, that, up to a renormalization of time, can be seen as random perturbations of an ordinary differential equation. 

\hspace{3mm}\newline
In Section 2 we will introduce our setting and give some examples of systems that fall into it.
Then, in Section 3 we will show that, under the hypothesis that the deterministic mean dynamic admits an interior attractor, the extinction time grows exponentially with the size of the system and that, when the system size goes to infinity, the limit set of the quasi-stationary distributions of the processes for the weak* convergence consists of invariant measures for the deterministic dynamic.
Finally in Section 4 we will study the support of these invariant limiting measures and prove that, under some additional large deviations hypotheses, their support lies within attractors bounded away from the extinction states. To do that we will compare two different notions of chain-recurrence, one given by the large deviations functional and the other a slight variation on Conley's $\delta$-T chain-recurrence.
Should the reader need some reminders about the objects used in this paper, he will find some basic properties and some references in the Appendix at the end.

\section{Model, notations and hypotheses}

\hspace{3mm}\newline

\noindent We denote by $\Delta$ the $d$-dimensional simplex.
\[ \Delta = \lbrace x \in \R^d \; ; \;\forall i=1 \cdots d \; \;  x_i \geqslant 0  \;\; \& \;\; \sum_{i=1}^d x_i =1 \rbrace \] 
We let $\iD$ denote the relative interior of $\Delta$ and
\[ \Delta_N=\Delta \cap \frac{1}{N} \Z^d\]  
\[\Delta^*_N=\iD \cap \frac{1}{N}\Z^d. \] Let $F:\Delta \to \R^d$ be a locally Lipschitz vector field such that :
\[ \forall x \in \Delta \quad \sum_{i =1}^d F_i(x)=0\]

Unless specified otherwise, the topology considered will be the topology induced by the classical $\R^d$ metric topology on $\Delta$.
Throughout the paper, if $A$ is a subset of a metric space $(E,d)$, we will denote by $N^\ep(A)$ its $\ep$-neighborhood
\[ N^\ep(A) =\lbrace x \in E \; ; \; d(x,A)<\ep \rbrace.\] 
We consider a family of Markov chains $(X_n^N)_{n \in \N}$ defined on a probability space $(\Omega, \mathcal{F}, \mathbb{P})$ taking values in the $d$-dimensional discrete simplex $\Delta_N$.

We denote by $\mathcal{F}_n^N$ the $\sigma$-algebra generated by $\lbrace X_i^N, i=1,...,n \rbrace$. For $A\in \mathcal{F}$ we let $\p_x[A]=\p[A \vert X_0=x ]$.

Throughout the paper the following hypothesis will always be assumed to hold.

\begin{sthyp}\label{base} 
The Markov process $X^N$ has the following properties :
\begin{enumerate}[(i)]
\item $\ds{X_{n+1}^N -X_n^N = \frac{1}{N} (F(X_n^N)+U_{n+1}^N)}$ \label{i}
\item $\mathbb{E}[ U_{n+1}^N \vert \mathcal{F}_n^N ]=0$ \label{ii}
\item There exists $\Gamma \geqslant 0$ such that $\Vert U_n^N \Vert \leqslant \sqrt{\Gamma}$ \label{iii}
\item  \label{iv} The boundary of the simplex is an \textit{absorbing} set:\begin{enumerate}[(a)]
\item for all $x \in \partial \Delta$ $\p_x [X^N_1 \in \partial \Delta ]=1$
 \item for all $x\in \Delta$ $\p_x[\exists n \, : \, X^N_n \in \partial \Delta ]=1$
 \end{enumerate}
\item  \label{v}$X^N$ restricted to $\Delta^*_N$ is \textit{irreducible}\[\forall x,y \in \Delta_N^*\quad \p_x[\exists n \, : \, X_n =y ]>0\]
and \textit{aperiodic} 
\[ \forall x \in \Delta^*_N \quad gcd( \lbrace n \; ; \; \p_x[ X_n =x ]>0 \rbrace )=1 \]

\end{enumerate}
\end{sthyp}

\begin{ex}[Guiding Thread]\label{fil}
Let $(p_{i,j}(x))_{i,j \in \lbrace 1 ... d\rbrace}$ be a family of real-valued continuous functions on $\Delta$  such that, for all $x \in \Delta$ :
\begin{gather}
 \forall i \neq j \quad p_{i,j}(x)=0 \Leftrightarrow x_i x_j=0, \label{a} \\
 p_{i,i}(x) = 0, \label{b}\\
  0\leqslant p_{i,j}(x) \leqslant 1,\label{c} \\
\sum_{i,j=1}^d p_{i,j}(x) \leqslant 1.\label{d}
\end{gather}

Let $(X_k^N)$ be the random walk on $\Delta_N$ defined by:
\[ \p \left[ X_{k+1}^N = X_{k}^N + \frac{1}{N}(e_j-e_i) \vert X^n_{k}=x \right]=p_{i,j}(x)\]
where $(e_i)_{i=1 \cdots d}$ is the canonical base of $\R^d$.
This type of model often occurs in population games. In this setting $N$ represents the size of the population. Each individual plays a pure strategy $i$ and $X^N$ represents then the vector of proportion of players of each strategy. The jump $ X_{k+1}^N = X_{k}^N + \frac{1}{N}(e_j-e_i)$ means that an individual switches his strategy from $i$ to $j$ at time $k$. The conditions on the family $(p_{i,j})$ mean that :
\begin{itemize}
\item At each time $k$, it is always possible that a player switches from his strategy $i$ to another strategy $j$ that is currently in use in the population.
\item No individual switches to an unused strategy. This makes sense for models based on strategy switching from imitations or models arising from ecology.
\end{itemize}

We define
\begin{align*} p_i(x) & =  \sum_{j=1}^d p_{j,i}(x)\\ q_i(x) & =  \sum_{j=1}^d p_{i,j}(x) \end{align*}
$p(x)$ the vector of coordinates $p_i(x)$ and $q(x)$ the vector of coordinates $q_i(x)$.

In this case $F(x) = p(x)-q(x)$. Hypothesis \ref{base}(\ref{i}) comes from the Markov property, \ref{base}(\ref{ii}) and \ref{base}(\ref{iii}) follow easily from the definition of the chain and \ref{base}(\ref{iv}) and \ref{base}(\ref{v}) follow from the fact that the functions $p_{i,j}$ are positive on the relative interior and vanish on the boundary.
\end{ex}

A class of examples that falls in Example \ref{fil} setting is given by \emph{Imitative Protocols games}, see e.g. \cite{San11}. Consider a population game with $d$ pure strategies, we let $U(x)=( U_1(x), \cdots U_d(x))$ denote the vector of payoffs when the population is in state $x\in \Delta$ and $\bar{U}(x)=\ds{\sum_{i=1}^d x_i U_i(x)}$ denotes the average payoff at population state $x$. Imitative protocols are of the form
\[ p_{i,j}(x)=x_i x_j r_{ij}(U(x),x) \]
with some additional assumptions on $r_{ij}$ to ensure that the problem is well-posed.
Under such a protocol, at each time we pick an individual uniformly at random among the population and give him a revision opportunity. The opportunity unfolds as such
\begin{itemize}
\item The individual picks an opponent uniformly at random among the population (he/she can pick him/her-self) and observes his/her strategy.
\item If the individual plays $i$ and the opponent plays $j$ then the individual switches from $i$ to $j$ with probability proportional to $r_{ij}$.
\end{itemize}

We now give some examples of Imitative Protocols

\begin{ex}[Pairwise Proportional Imitation]

After selecting an opponent the agent imitates only if the opponent's payoff is higher than his own, doing so with probability proportional to the payoff difference.
\[ p_{i,j}(x)= x_i x_j (U_j(x)-U_i(x))^+ \]
where $(y )^+$ stands for $\max (y , 0)$. 

The mean dynamic generated by this protocol is
\begin{align*}
\dot{x} & =  \sum_{j=1}^d  p_{j,i}(x) -  \sum_{j=1}^d p_{i,j}(x) \\
        & = \sum_{j=1}^d x_j x_i \left( U_i(x)-U_j(x) \right)^+ - x_i \sum_{j=1}^d \left( U_j(x) -U_i(x) \right)^+ \\
        & = x_i \sum_{j=1}^d x_j (U_i(x) -U_j(x)) \\
        & = x_i(U_i(x) - \bar{U}(x))
\end{align*}
        
We then get the well known replicator dynamic, a dynamic extensively studied in ecology and evolutionary game theory, see e.g. \cite{HofSig98}
\end{ex}

\begin{ex}[Aspiration and Random Imitation]

A particular case of the former example is the aspiration and random imitation model, see e.g. \cite{BenWei03}, \cite{BinGalSam95}, \cite{BjoWei96} and \cite{BinSam97}. At each time we pick an individual at random in the population and look at his/her "satisfaction", a payoff-like function $u_i(x)$ where $i$ is the type of the drawn individual. If this satisfaction is lower than a certain aspiration level then switch to another type chosen at random in the population, otherwise stay at current type. The morality of this model is that, if your type isn't performant enough for your tastes then switch to another type. The aspiration levels are independent random variables uniformly distributed on intervals $[a_i(x),b_i(x)]$ with $a_i(x) \leqslant u_i(x) \leqslant b_i(x)$. This model gives us 
\[ p_{i,j}(x)=x_i x_j \frac{b_i(x)-u_i(x)}{b_i(x)-a_i(x)} \]
If we assume that the aspiration level bounds are not type-dependant, meaning $a_i(x)=a(x)$ and $b_i(x)=b(x)$ for all $i$ then we get a mean field given by
\[ \dot{x}_i =x_i \frac{u_i(x) - \sum_{j=1}^d x_j u_j(x)  }{b(x)-a(x)} \]
which is a replicator dynamic with fitness functions $f_i(x)=\frac{u_i(x)}{b(x)-a(x)}$

Alternatively we can assume that the aspiration levels follow the type payoff by the relation $b_i(x)=\beta_i u_i(x)$ and $a_i(x)=\alpha_i u_i(x)$ with $\alpha_i < 1 < \beta_i$. In this case we get a dynamic
\[\dot{x}_i=x_i \left( \sum_{j=1}^d x_j v_j - v_i \right)\]
where $v_j=\frac{\beta_j-1}{\beta_j-\alpha_j}$.

Again this is a replicator dynamics with fitness function $f_i(x)=-v_i$.
\end{ex}

\begin{ex}[Imitation Driven by Dissatisfaction]

In this protocol, when a $i$ player receives a revision opportunity, he opts to switch strategies with a probability that that is linearly decreasing in his current payoff. Should he decide to change, then he will imitate a randomly selected opponent. This protocol also falls under the former example of aspiration games with constant aspiration levels $a_i(x)=A$, $b_i(x)=B$.
This gives the following dynamic
\[ \dot{x}_i = \frac{1}{B-A}(F_i(x) - \bar{F}(x)) \]
Again we recognize a replicator dynamic.
\end{ex}

For other examples see e.g. \cite{San11}.

\section{Convergence of QSD and absorption time}

\hspace{1mm} \newline

We denote by $\lbrace\varphi_t\rbrace$ the flow induced by $F$. In order to compare the trajectory of $\varphi_t$ with those of $(X_n^N)$ it's convenient to introduce the continuous process $\hat{X}^N : \R \to \R^m$ defined by \[\hat{X}^N(k/N)=X_k^N \quad \forall k \in \N\]
and extended on every interval $[k/N,(k+1)/N]$ by piecewise linear interpolation.

Let
\[D^N(T)=\max_{0\leqslant t \leqslant T} \Vert \hat{X}^N(t)-\varphi_t(X_0^N)\Vert\]
be the variable measuring the distance between the trajectories $t \mapsto \hat{X}^N(t)$ and $t \mapsto \varphi_t(X_0^N)$.

We recall this convergence theorem of Benaïm and Weibull.
\begin{thm}\label{champ}
For every $T >0$, there exists $c >0$ (depending only on $F$,$\Gamma$ and $T$) such that, for every $\ep >0$, and for $N$ large enough :
\[\mathbb{P} [D^N(T) \geqslant \varepsilon ] \leqslant 2d\e{-\ep^2cN}\]
\end{thm}

For a detailed proof of this result, see \cite{BenWei03}.

\vspace{3mm}


We define $T^N_0$ to be the \textit{absorption time}. 
\[T^N_0=\inf \lbrace n>0 ; X^N_n \in \partial \Delta \rbrace\]
Hypothesis \ref{base}(\ref{iv}) implies that, whatever the initial state is, the process will almost surely be absorbed, i.e.
\[ \forall x \in \Delta^N \qquad \p_x [T^N_0 < \infty ]=1.\]

A probability measure $\mu$ on the discrete relative interior of the simplex $\Delta_N^{*}$ is said to be a \textit{quasi-stationary distribution}, thereafter referred as QSD, if and only if, for every Borelian set $A \subset \Delta_N^{*}$ and every $n>0$,
\[ \p_\mu [X_n^N \in A \vert T^N_0 > n ] = \mu(A).\]
We remark that, in this case, $\mu$ is a fixed point for the conditional evolution 
\[\nu \mapsto \p_\nu [ X_n^N \in \centerdot \vert T_0 > n ]\]

The following proposition is a classic QSD result and follows easily from the Perron-Frobenius theorem. 
\begin{prop}\label{abso}
For every $N$, there exists an unique quasi-stationary distribution $\mu^N$ obtained as the only left eigenvector $\mu^N$ of the transition matrix of the Markov chain restricted to $\Delta_N^{*}$ verifying
\[\forall i \in \lbrace 1, \cdots d \rbrace \quad \mu^N_i >0 \quad ; \quad \sum_i \mu^N_i =1\]
The corresponding eigenvalue $0 <\rho_N= \e{-\theta_N}<1$ is such that
\[\p_{\mu^N} [T^N_0 > n ] = \e{- \theta_N n} \] 
Hence, starting from $\mu^N$, the expectation of $T_0^N$ is 
\[ \E_{\mu^N} [T^N_0] = \frac{1}{1-\rho_N} \]
\end{prop}

\begin{pre}

For a detailed proof see e.g. \cite{MelVil11}.

\end{pre}

\subsection{Absorption time}

\hspace{3mm}\newline

\noindent A set $A\subset \Delta$ is called an \textit{attractor} for the flow $\lbrace\varphi_t\rbrace$ if
\begin{enumerate}[(i)]
\item $A$ is compact and invariant, i.e. for every $t\in \R$  $\varphi_t(A)=A$.
\item There exists a neighborhood $U$ of $A$, called a fundamental neighborhood, such that
\[ \lim_{t \to \infty} d(\varphi_t(x),A)=0 \]
uniformly in $x$ in $U$.
\end{enumerate}

\begin{thm} \label{attra}
Starting from $\mu^N$, the law of the absorption time and its expectation are given by Proposition \ref{abso}. If we further assume that the flow $\lbrace \varphi_t \rbrace$ admits an attractor $A \subset \iD$, then, there exists $\gamma>0$ such that the following estimate holds :
\[  0\leqslant 1-\rho_N \leqslant  O\left( \frac{ \e{-\gamma N}}{N} \right)\]
Thus, there exists a constant $C>0$ such that 
\[\E_{\mu^N} [ T_0 ] \geqslant C N \e{\gamma N} \]
\end{thm}

\begin{pre}

Let $V \subset \mathring{\Delta}$ and let $k\in\N$. By the QSD property we have:
\begin{align*} \rho_N^k \mu^N(V) & = \sum_{x \in \Delta_N }p_N^k(x,V)\mu^N(x) \\
& \geqslant \sum_{x \in V\cap\Delta_N }p_N^k(x,V)\mu^N(x)\\ & \geqslant   \inf_{x \in V\cap\Delta_N} p_N^k(x,V) \mu^N(V). \end{align*}

Thus \[\rho_N^k \geqslant \inf_{x \in V\cap\Delta_N} p_N^k(x,V).\]

Let $U\subset \iD$ be a fundamental neighborhood of the attractor $A$. We know that $d(\varphi_t(x),A)$ converges uniformly to $0$ over $U$. Hence 
\[\forall \varepsilon >0 \quad \exists T(\varepsilon)>0 \quad \forall t \geqslant T(\varepsilon)\quad \forall x \in U \quad d(\varphi_t(x),A)<\varepsilon .\]

Let $\alpha=d(A,U^c)$, $\varepsilon < \alpha$, $T=T(\varepsilon)$ and $\delta < \alpha-\varepsilon$.

For all $x \in U \cap\Delta_N$
\begin{align*} p_N^{[NT]}(x,U^c) & \leqslant  \p_{x} [ d(X^N_{[NT]},A) >\alpha ]\\ & \leqslant   \p_{x} [ d(X^N_{[NT]},\varphi_T(x)) >\alpha-\varepsilon ]\\ & \leqslant  \p_{x} [ d(X^N_{[NT]},\hat{X}^N_T) +  d(\hat{X}_T,\varphi_T(x)) >\alpha-\varepsilon ]\\ & \leqslant   \p_x \left[ D_N(T)>\alpha -\varepsilon - ( \Vert F \Vert + \sqrt{\Gamma}) \frac{NT-[NT]}{N} \right] \\ & \leqslant  2d\e{-\delta^2 c N} \text{ for }N\text{ large enough (see Theorem \ref{champ})} \end{align*}

Then
\begin{align*} \rho_N^{[NT]}  & \geqslant \inf_{x \in U \cap \Delta_N} p_N^{[NT]}(x,U) \\
& \geqslant  1- \max_{x \in U\cap \Delta_N} p_N^{[NT]}(x,U^c)\\  & \geqslant  1- 2d\e{-\delta^2 c N} \end{align*}

Therefore \[1-\rho_N \leqslant 1-\left(1- 2d\e{-\delta^2 c N}\right)^{\frac{1}{[NT]}}\]

\begin{align*}\left(1- 2d\e{-\delta^2 c N}\right)^{\frac{1}{[NT]}}  & =  \e{{\frac{1}{[NT]}}\log{\left(1- 2d\e{-\delta^2 c N}\right)}} \\ & =   \e{{\frac{-1}{[NT]}}2d\e{-\delta^2 c N}}+o(\e{-\delta^2 c N}))\\ & =  1- \frac{2d\e{-\delta^2 c N}}{[NT]} + o\left(\frac{2d\e{-\delta^2 c N}}{[NT]}\right) \end{align*}

In conclusion we have 
\[  0\leqslant 1-\rho_N \leqslant  O\left( \frac{ \e{-\gamma N}}{N} \right)\]

\end{pre}

\subsection{Convergence of the QSD to an invariant measure}

\hspace{3mm}\newline

\noindent A probability measure $\mu$ on $\Delta$ is called an \textit{invariant measure} for the flow $\lbrace \varphi_t \rbrace$ if, for all $t \in \R$ and all borelian set $A \in \mathcal{B}(\Delta)$, $\mu(\varphi_t^{-1}(A))=\mu(A)$. 

\begin{thm}\label{liminvar}
We suppose that the flow $\lbrace \varphi_t \rbrace$ admits an attractor $A \subset \mathring{\Delta}$. 
Then the set of limit points of $\lbrace \mu^N \rbrace $  for the weak* topology is a subset of the set of invariant measures for the flow $\lbrace\varphi_t\rbrace$.
\end{thm}

\begin{pre}

Let $f$ be a Lipschitz function from $\Delta$ to $\R$ with constant $L$. We suppose that the sequence $\mu^N$ weakly converges to a measure $\mu$. Let $t>0$.
We want to prove that \[\displaystyle{\lim_{N \to \infty} \int f(x) \mu^N(dx) - \int f(\varphi_t(x))\mu^N(dx)=0}\]
The QSD property gives us that, for all $k$ 
\[\E_\mu [f(X_n)]= \int f(x) \mu^N(dx) =  \int \E_x \left[ f(X_{k}^N) \bigg\vert T^N_0 > k \right]\mu^N(dx) \]

Let
\[I  =  \left| \int f(x) \mu^N(dx) - \int f(\varphi_t(x))\mu^N(dx) \right| \]
Then, for all $k$,
\begin{align*}
I & =  \left| \int f(x) \mu^N(dx) - \int f(\varphi_t(x))\mu^N(dx) \right| \\
 & =  \left| \int \E_x \left[ f(X_{k}^N) \bigg\vert T^N_0 > k \right]  \mu^N(dx) - \int f(\varphi_t(x))\mu^N(dx) \right|  \\
  & =  \left| \int \E_x \left[ f(X_{k}^N) - f(\varphi_t(x)) \bigg\vert T^N_0 > k \right] \mu^N(dx) \right|  \end{align*}
  
In particular, for $k=[Nt]$.
  
\begin{align*}   I & =  \left| \int \E_x \left[ f(X_{[Nt]}^N)  - f(\varphi_t(x)) \bigg\vert T^N_0 >[Nt] \right] \mu^N(dx) \right|
\end{align*}

By Theorem \ref{champ}, we know that, for $N$ large enough, we have  
\[ \p_x [ D_N(t)>\delta ] \leqslant 2d\e{-\delta^2 c N}. \]

Thus \[\E_x [D_N(t) ] = \int_0^{+\infty}   \p_x [ D_N(t)>\delta ] d\delta \leqslant  \int_0^{+\infty} 2d\e{-\delta^2 c N} d\delta = \frac{d\sqrt{\pi}}{\sqrt{cN}}\]

From Theorem \ref{attra}, we can infer $1-\rho_N \leqslant O\left( \frac{ \e{\gamma N}}{N} \right)$.

This implies $\e{\theta_N [Nt]} -1 \underset{N \to +\infty}{\longrightarrow} 0$, which gives us the boundedness of $\e{\theta_N [Nt]}$.

Hence
\begin{align*} I & =  \left| \int \E_x \left[ f(X_{[Nt]}^N)  - f(\varphi_t(x)) \bigg\vert T^N_0 > [Nt] \right] \mu^N(dx) \right| \\
 & \leqslant  \left| \int \frac{\E_x \left[ f(X_{[Nt]}^N)  - f(\varphi_t(x)) \right]}{\p_x \left[ T_0^N >[Nt] \right]} \mu^N(dx) \right| \\
 & \leqslant  \left| \int \frac{\E_x \left[ L \vert X_{[Nt]}^N -\hat{X}_t +\hat{X}_t  - \varphi_t(x) \vert \right]}{\p_x \left[ T_0^N >[Nt] \right]} \mu^N(dx) \right| \\
 & \leqslant  \left| \int \frac{\E_x \left[ L (D_N(t) + (Nt -[Nt])\frac{1}{N}(\Vert F \Vert + \Gamma)) \right]}{\p_x \left[ T_0^N >[Nt] \right]} \mu^N(dx) \right| \\
 & \leqslant  \left| L \left( \frac{d \sqrt{\pi}}{\sqrt{cN}}+ \frac{\Vert F \Vert + \Gamma}{N} \right) \e{\theta_N [Nt]} \right| \underset{N \to +\infty}{\longrightarrow} 0
\end{align*}

\end{pre}

\section{Support of the limiting measure}

\hspace{3mm}\newline

\noindent Let $L=L(\lbrace\mu_N \rbrace)$ denote the limit set of the sequence $(\mu_N)_{N \in \N}$ for the weak* topology. In view of Theorem \ref{liminvar}, $L$ consists of invariant measures. As the QSD have their support inside $\iD$, it is natural to study whether the limiting measure also take their support in $\iD$.  However, by the Poincaré Recurrence Theorem, every $\mu \in L$ is supported by the Birkhoff center
\[BC(\varphi)=\adh{\lbrace x \in \Delta \; ; \; x \in \omega(x) \rbrace} \] 
Since the boundary of $\Delta$ intersects the Birkhoff center (e.g the vertices of the simplex are equilibria and thus inside the Birkhoff center), knowing that the QSD converges to an invariant measure is not enough, we have to further study the support of the measure $\mu$ to ensure that it is strictly inside the interior of the simplex.
For that we will need large deviation assumptions.

\begin{hyp}\label{H1'}
For all $\alpha>0$, there exists a function 
\begin{align*}
S : V_\alpha \times \R \times  \mathcal{C}_x([0,T],V_\alpha) & \longrightarrow  \bar{\R}_+ \\
 (x,T,\phi) & \longmapsto  S(x,T,\phi)
\end{align*}
with the following properties, where $V_\alpha= \Delta \setminus \adh{N^\alpha (\partial \Delta)}$, $\mathcal{C}([0,T],V_\alpha)$ is the set of continuous functions $\psi$ from $[0,T]$ to $V_\alpha$ and $\mathcal{C}_x([0,T],V_\alpha)$ is the set of continuous functions $\psi$ from $[0,T]$ to $V_\alpha$ such that $\psi(0)=x$, both equipped with the topology of uniform convergence.

\begin{itemize}
\item For every $s \in ]0,\infty[$ and $T>0$, the set \[\lbrace \phi \in  \mathcal{C}_x([0,T],V_\alpha) \text{ s.t. } S(x,T,\phi) \leqslant s \rbrace\] is a compact set
\item For $x\in \iD$ and $T>0$, $S(x,T,\phi)=0$ $\Leftrightarrow$ $\dot{\phi}_s=F(\phi_s)$ $\forall s \in [0,T]$
\item $\hat{X}$ satisfies a large deviation principle with rate function $S$ and speed $1/N$ uniformly in $x$ on compact subsets of $V_\alpha$,

 i.e. for $K$ compact subset of $V_\alpha$, $T>0$ and $A \subset  \mathcal{C}_x([0,T],V_\alpha)$
\begin{align*} - \sup_{x \in K} \inf_{\phi \in \mathring{A}} S(x,T,\phi) & \leqslant \liminf_{N \to \infty} \frac1N log \inf_{x \in K} \p_x [ \hat{X}^N \in A ] \\ \limsup_{N \to \infty} \frac1N log \sup_{x \in K} \p_x [ \hat{X}^N \in A ]& \leqslant - \inf_{x \in K} \inf_{\phi \in \bar{A}} S(x,T,\phi)\end{align*}
In particular when $K=\lbrace x \rbrace$ we get the "classical" large deviation principle\begin{align*} -  \inf_{\phi \in \mathring{A}} S(x,T,\phi) & \leqslant \liminf_{N \to \infty} \frac1N log  \p_x [ \hat{X}^N \in A ] \\ \limsup_{N \to \infty} \frac1N log  \p_x [ \hat{X}^N \in A ]& \leqslant -  \inf_{\phi \in \bar{A}} S(x,T,\phi)\end{align*}
\item $S$ is linear with regards to the concatenation of functions, i.e. if $\tilde{T}<T$
\[ S(x,T, \phi)= S(x,\tilde{T},\restriction{\phi}{[0,\tilde{T}]}) + S(\phi(\tilde{T}),T-\tilde{T}, \restriction{\phi}{[\tilde{T},T-\tilde{T}]}) \]
\item $\ds{\lim_{T \to 0} S(x,T,\phi)=0}$ uniformly in $x\in V_\alpha$ and $\phi \in \mathcal{C}_x([0,T],V_\alpha)$.
\end{itemize}
\end{hyp}

\begin{hyp}\label{H2'}
$\forall c >0$ $\exists U_{c}$ an open neighborhood of $\partial \Delta$ such that
\[\lim_{N \to \infty} \inf_{x \in U_c} \frac1N log \p_x [ \hat{X}^N(1) \in \partial \Delta ] \geqslant -c\]
\end{hyp}

\begin{ex}
Before going further we will verify that Hypotheses \ref{H1'} and \ref{H2'} hold for the nearest neighbor random walk model introduced in Example \ref{fil}.
%
%
%
%

To show that Hypothesis \ref{H1'} holds we will use Theorem 6.3.3 in \cite{pDUP97a} with Conditions 6.2.1 and 6.3.1, which will gives us a Laplace principle for $\restriction{\hat{X}^N}{[0,T]}$ that holds uniformly on compacts, and Theorem 1.2.3 in \cite{pDUP97a} will, in turn, give us the desired uniform on compacts large deviation principle.

We only need to show that Conditions 6.2.1 and 6.3.1 holds

When transcribing our model in \cite{pDUP97a} setting, we get, for $x \in V_\alpha$
\[ \mu(dy \vert x ) = \sum_{i\neq j} p_{i,j}(x) \delta_{e_j-e_i}(dy) + \left(1- \sum_{i\neq j} p_{i,j}(x)\right) \delta_0(dy) \]

We define
\[ H_\mu(x,\alpha) = log \int_{\R^d} exp\langle\alpha,y\rangle \mu(dy\vert x) \]
and 
\[ L_\mu(x, \beta) = \sup_{\alpha \in \R^d} \lbrace \langle \alpha,\beta \rangle - H_\mu(x, \alpha) \rbrace \]
We refer to Chapter 6.2 in \cite{pDUP97a} for elementary properties of these functions.

\begin{defin}
$\mu(dy\vert x)$ is said to verify Condition A (called 6.2.1 in \cite{pDUP97a})  if 
\begin{enumerate}[(i)]
\item For each $\alpha \in \R^d$, $\ds{\sup_{x \in \R^d} H_\mu(x,\alpha) < \infty}$
\item The function mapping $x \in \R^d \mapsto \mu(\centerdot \vert x )$ is continuous in the topology of weak convergence.
\end{enumerate}
\end{defin}

\begin{defin}
$\mu(dy\vert x)$ is said to verify Condition B ( called 6.3.1 in \cite{pDUP97a}) if 
\begin{enumerate}[(i)]
\item The relative interior of the convex hull of the support of $\mu(dy\vert x)$, $Ri(Conv(Supp(\mu(\centerdot\vert x))))$, doesn't depend on $x$.
\item $0 \in Ri(Conv(Supp(\mu(\centerdot\vert x))))$
\end{enumerate}
\end{defin}

Let $ \alpha >0$, we'll now prove the L.D.P. on $V_\alpha$
The conditions, as they are written, demand that $x$ may take values in all of $\R^d$ while, in our model, we only use $\Delta$. To remedy to that we'll first embed our $d$-dimensional simplex in $\R^{d-1}$ and then extend $\mu$ to a kernel $\eta$ defined on all of $\R^{d-1}$ by taking $\eta(dy\vert x) = \mu(dy \vert p_\alpha(x))$, where $p_\alpha$ is the convex projection on $\adh{V_\alpha}$.

This way we have a probability kernel $\eta$ that is defined on all of $\R^{d-1}$, it is then easy to verify that Conditions A and B hold for $\eta$. We thus get a LDP for $\eta$ on all of $\R^{d-1}$ with speed $1/N$ and rate function 
\[S_\eta= \begin{cases} \int_0^T L_\eta(\phi(t)  ,\dot{\phi}(t)) dt \text{ if }\phi\text{ is uniformly continuous}\\
0 \text{ otherwise } \end{cases} \]

We now remark that, when $x \in V_\alpha$, $\eta(dy\vert x) =\mu(dy\vert x)$ and thus $H_\eta(x,v)=H_\mu(x,v)$ and $L_\eta(x,u)=L_\mu(x,u)$.

From that we deduce that, if $\phi \in \mathcal{C}([0,T],V_\alpha)$, then $L_\eta(\phi,\dot{\phi})=L_\mu(\phi,\dot{\phi})$ and finally $S_\eta(x,T,\phi)=S_\mu(x,T\phi)$.

We finally get that Hypothesis \ref{H1'} holds for our nearest neighbor random walk with the following rate function.
\[ S(x,T,\phi)=\begin{cases} \int_0^T L_\mu(\phi(t)  , \dot{\phi}(t)) dt \text{ if }\phi\text{ is uniformly continuous}\\
0 \text{ otherwise } \end{cases} \]
We still have to verify that the wanted properties holds for this rate function: 
\begin{itemize}

\item From Proposition 6.2.4 in \cite{pDUP97a} we get that, for every $s \in ]0,\infty[$ and $T>0$, the set \[\lbrace \phi \in  \mathcal{C}_x([0,T],\Delta) \text{ s.t. } S(x,T,\phi) \leqslant s \rbrace\] is a compact set.
\item Let $x\in V_\alpha$, $T>0$, it is already known that
\[S(x,T,\phi)=0 \quad \Leftrightarrow \quad \dot{\phi}_s=F(\phi_s) \qquad\forall s \in [0,T]\]
\item The LDP comes from Theorem 6.3.3 in \cite{pDUP97a}
\item The linearity of $S(x,T,\phi)$ follows easily from it's definition as an integral.\end{itemize}

Let's now prove that Hypothesis \ref{H2'} holds too.

\begin{prop}
Suppose that, for every couple $(i,j) \in \lbrace 1 , \cdots d \rbrace^2$ there exists $k \in \N$ such that $\ds \frac{\partial^k p_{i_j}(x)}{\partial x_i^k} \neq 0$ on $\lbrace x \in \Delta \; ; \; x_i=0 \rbrace$. Then  Hypothesis \ref{H2'} holds
\end{prop}

\begin{pre}

We want to ensure that, for all $c>0$ there exists an open neighborhood $U_c$ of $\partial \Delta$ such that
\[\lim_{N \to \infty} \inf_{x \in U_c} \frac1N log \p_x [ \hat{X}^N(1) \in \partial \Delta ] \geqslant -c.\]

We know that $q_i(x)$ and $p_i(x)$ go to $0$ as $x_i$ goes to $0$ and that there exists $k \in \N$ such that $\ds \frac{\partial^k q_i(x)}{\partial x_i^k} > 0$ on $\lbrace x \in \Delta \; ; \; x_i=0 \rbrace$ and thus $\ds \frac{\partial^k q_i(x)}{\partial x_i^k} > a$ on a sufficient small neighborhood of $\lbrace x \in \Delta \; ; \; x_i=0 \rbrace$ with $a>0$.

Let $1>b>0$ and let $ x\in U_c= \lbrace x \; ; \; \exists i \; x_i<b\rbrace$. We have 
\begin{align*} \p_x [ \hat{X}^N(1) \in \partial \Delta ] & \geqslant  \p_x [ X^N_{[Nb]} \in \partial \Delta ] \\ & \geqslant  \ds{\prod_{j=1}^{[Nb]} q_i(x^{(j)})}\end{align*}
where $(x^{(j)})_j$ is a sequence of points in $\Delta_N$ such that $x_i^{(j)}=\frac{j}{N}$. 

If $b$ is small enough we have, for $j \leqslant [Nb]$ 
\[q_i(x^{(j)})\geqslant a \left(\frac{j}{N}\right)^k\]

Then 
\[\prod_{j=1}^{[Nb]}q_i(x^{(j)})  \geqslant a^{[Nb]}\frac{[Nb]!}{N^{[Nb]}} \]

Thus 
\begin{align*}\inf_{x \in U_c} \frac1N log \p_x [ \hat{X}^N(1) \in \partial \Delta ]  \geqslant & C\left( \frac{[Nb]}{N}(log(a)-1) \right. \\ &  \left. \qquad + \frac{[Nb]}{N}log{\frac{[Nb]}{N} + \frac{1}{2N}log(2 \pi [Nb])}\right)
\end{align*}
where $C$ is a constant

As $N$ goes to infinity, the right-hand term goes to $Cb(log(a)-1) +Cblog(b)$ which is greater than $-c$ for $b$ small enough.
Thus, for $N$ large enough
\[\lim_{N \to \infty} \inf_{x \in U_c} \frac1N log \p_x [ \hat{X}^N(1) \in \partial \Delta ] \geqslant -c\]

Hence our random walk model satisfies Assumption \ref{H2'}.
\end{pre}

Finally both hypotheses holds for our model.

\end{ex}

\begin{defin}
We define 
\[\L (x,y)= \ds{\limsup_{T \to \infty} \inf_{\phi \in \mathcal{C}_x([0,T]), \, \phi(T)=y} S(x,T,\phi)}\]
We will say that $x$ $\L$-leads to $y$ (denoted by $x\rightsquigarrow_\L y$) if, for every $\ep>0$, there exists a path of points $x=\xi_0, \xi_1, \cdots \xi_{n(\xi)}=y$ such that 
\[A_{n(\xi)}(\xi)=\ds{ \sum_{k=1}^{n(\xi)} \L(\xi_k,\xi_{k+1})} < \ep\]

We define $\ds{B_\L (x,y) =\inf_{\xi \text{ linking }x\text{ to }y} A_{n(\xi)}(\xi)}$.

Thus $x\rightsquigarrow_\L y$ iff $B_\L(x,y)=0$

We will say that $x$ is $\L$-chain recurrent if $x \leadsto_{\L} x$ and will denote by $\mathcal{R}_{\L}(\varphi)$ the set of all $\L$-chain recurrent points.

If $x$ and $y$ are two points of $\mathcal{R}_{\L}$ verifying $x \leadsto_{\L} y$ and $y \leadsto_{\L} x$ we will then denote $x \sim_{\L} y$. The equivalence classes for this relation will be called $\L$-basic classes. We define a partial order on these classes by $[x] \prec_{\L} [y]$ if $x \leadsto_{\L} y$. A maximal $\L$-basic class will be called a $\L$-quasi-attractor.
\end{defin}

\begin{hyp}\label{H3'}
There is only a finite number of $\L$-basic classes in $\mathring{\Delta}$ denoted by $K_i \, , \, i =1\cdots \nu$. We suppose that they are closed sets and indexed in such a way that the $k$ first $\lbrace K_i \rbrace_{i =1 \cdots k}$ are the $\L$-quasi-attractors  and the $\nu-k$ others aren't.
\end{hyp}

\begin{prop}\label{cont}
The function $\L$ has the following properties :
for every sequence $(x_n)_{n \in \N} \in V_\alpha$ converging to an $x\in V_\alpha$ and every $y\in V_\alpha$ we have
\[\lim_{n\to \infty } \L(x_n ,y) =\L(x,y) \qquad \lim_{n\to \infty } \L(y,x_n ) =\L(y,x)\]
\end{prop}

\begin{pre}

This proposition follows easily from Hypothesis \ref{H1'}

\end{pre}

\vspace{3mm}

The following theorem is the main result of this section, giving us more insight in the support of the limiting measure $\mu$.

\begin{thm}\label{main}
We suppose that the flow $\lbrace \varphi_t\rbrace$ associated with the mean dynamic $\dot{x}=F(x)$ has an interior attractor. Under Hypotheses \ref{H1'}, \ref{H2'} and \ref{H3'}, we have :

The limiting measure $\mu$ has its support in the union of the $\L$-quasi-attractor, i.e. in $\ds{ \bigcup_{i=1}^k K_i }$
\end{thm}

We will prove this theorem under an intermediary set of hypotheses then we will prove that the announced hypotheses imply the intermediary hypotheses.

\subsection{Absorption-preserving pseudo-orbit}
\hspace{1mm}\newline
Here we introduce a different notion of chains for our dynamical system using absorption preserving $\delta, T$ pseudo-orbit, an analog to $\delta,T$ pseudo-orbits introduced by Conley \cite{Con78}, which have been extensively studied in the past.

\begin{defin}
Let $\lbrace\varphi_t\rbrace_{ t \in \R}$ be a flow given by an ordinary differential equation on $(\Delta,d)$ for which $\partial \Delta$ is an invariant set.
We will call $(\delta ,T)$ absorption preserving pseudo-orbit ($\delta,T$-ap-pseudo orbit) from $x$ to $y$ a piecewise continuous path 
\begin{align*} x=x_0,  & \lbrace\varphi_t(x_1)\, ;\, t \in [0,t_1] \rbrace, \lbrace\varphi_t(x_2)\,;\, t \in [0,t_2] \rbrace \\ & \cdots \lbrace\varphi_t(x_1)\,;\, t \in [0,t_{k}] \rbrace, x_{k+1}=y \quad k\geqslant 1\end{align*}
which is uniquely defined by the sequences of points $x_0, \cdots x_{k+1}$ and times $t_1, \cdots t_{k}$ such that the following hypotheses hold:
\[\begin{cases} d(x,x_1)<\delta \\
d(\varphi_{t_j}(x_j),x_{j+1}) <\delta \quad \forall j=1 \cdots k \\
t_i \geqslant T \quad \forall i=1\cdots k \\
x_j \in \partial \Delta \Rightarrow x_{j+1} \in \partial \Delta \quad \forall j=0 \cdots k \\
\end{cases} \]
We will denote then $x \leadsto_{ap,\delta,T}  y$

If $x \leadsto_{ap,\delta, T} y$ for every $\delta >0 $ and every $T>0$, we will denote $x \leadsto_{ap} y$

The point $x$ will be said to be ap-chain recurrent if $x \leadsto_{ap} x$, we define $\mathcal{R}_{ap}(\varphi)$ the set of ap-chain recurrent points.

If $x$ and $y$ are two points of $\mathcal{R}_{ap}(\varphi)$ such that $x \ap y$ and $y \ap x$  we will write $x \sim_{ap} y$. The equivalence classes for this relation will be called ap-basic classes. We define a partial order on these classes by $[x] \prec_{ap} [y]$ if $x \ap y$. A maximal ap-basic class will be called a ap-quasi-attractor.
\end{defin}

The ap-chain recurrent points and the ap-basic classes have some interesting properties which will be of use when proving the result on the support of the limit measure $\mu$. We enumerate and prove some of them.

\begin{prop}
Let $x \in \mathcal{R}_{ap}$. Then $\bar{[x]}_{ap} \subset [x]_{ap} \cup \partial \Delta$.
Moreover, for all $t>0$ $[x]_{ap}$ is $\varphi_t$ invariant. 
\end{prop}

\begin{pre}

Let $y \in \bar{[x]}_{ap}$ and let $(y_k)_{k \in \N}$ be a sequence of elements of $[x]_{ap}$ converging to $y$. Suppose $y \not\in \partial \Delta$.

Let $\delta >0$ and $T>0$ and let $k$ such that $d(y_k,y)<\delta$.

There exists a $\delta ,T$ pseudo-orbit $x, (x_1,t_1) , \cdots , (x_n, t_n) ,y_k $ linking $x$ to $y_k$. Thus $x, (x_1,t_1) , \cdots , (x_n, t_n) ,y $ is a $2\delta ,T$ pseudo-orbit linking $x$ à $y$.

Hence $x \ap y$. 

A similar reasoning on $\delta,T$-pseudo-orbit linking $y_k$ to $x$ gives us $y \ap x$ and then $y \in [x]_{ap}$.

We show now that $[x]_{ap}$ is an invariant set.
 
Let $T,T',\ep>0$. $\varphi_{T'}$ is an uniformly continuous application. Let then $\delta <\ep$ such that $d(x,y)<\delta$ implies  $d(\varphi_{T'}(x),\varphi_{T'}(y)) <\ep$.

We know that there exists a $\delta, T$ pseudo-orbit  $x, (x_1,t_1) , \cdots , (x_n, t_n) ,x $ linking $x$ to $x$. Hence $x, (x_1,t_1) , \cdots , (x_n, t_n+T') ,\varphi_{T'}(x) $ is a $\ep, T$ pseudo-orbit linking $x$ to $\varphi_{T'}(x)$. 

Before proving the converse, let's just remark that, if $\delta_1 < \delta_2$ and $T_1 > T_2$ then every  $\delta_1, T_1$ pseudo-orbit is also a $\delta_2, T_2$ pseudo-orbit.

Let's now suppose that $T > T'$ 

We consider again our $\delta, T$ pseudo-orbit $x, (x_1,t_1) , \cdots , (x_n, t_n) ,x $ linking $x$ to $x$.
One can remark that $t_1+t_2-T' > T$.

Thus $\varphi_{T'}(x), (\varphi_{T'}(x_1),t_1+t_2-T), (x_3,t_3) , \cdots , (x_n, t_n) ,x $ is a $\ep ,T$ pseudo-orbit linking $\varphi_{T'}(x)$ to $x$ if $\delta$ is chosen small enough.

Then $x \ap \varphi_{T'}(x)$.

By composing the $\delta ,T$ pseudo-orbits linking $x$ and $\varphi_{T'}(x)$ by $\varphi_{-T'}$ we get that $x \sim_{ap} \varphi_T(x)$ for every $T$ in $\R$.

\end{pre}

\begin{prop}
Let $x \in \Delta$. If $x \in \partial \Delta$ or $\omega (x) \subset \iD$ then $\omega(x) \subset \mathcal{R}_{ap}$. From this we get that, for every $x$ in $\Delta$, $\omega (x) \cap \mathcal{R}_{ap} \not= \emptyset$.
\end{prop}

\begin{pre}

The first point is a well-known result for "classic" chain-recurrence and easily extended to ap-chain-recurrence. Furthermore, if $x \in \iD$ and $\omega(x) \cap \partial \Delta \not= \emptyset$ then, by taking $y \in \omega(x) \cap \partial \Delta$, we get $\omega(y) \subset \omega(x)$ and $\omega(y) \subset \mathcal{R}_{ap}$.

\end{pre}

\begin{prop}
If $[x]_{ap}$ is maximal, then $x \ap z$ if and only if $z \in [x]$. As a consequence we also get that every quasi-attractor is a closed set.
\end{prop}

\begin{pre}

This result is trivial as soon as $z \in \mathcal{R}_{ap}$, thus we only have to prove it for $z \not \in \mathcal{R}_{ap}$. We know that $\omega (z) \cap \mathcal{R}_{ap} \not= \emptyset$ and that, if $u \in \omega(z)\cap \mathcal{R}_{ap}$ then $z \ap u$ thus $x \ap u$. Hence $u \in [x]_{ap}$, from that we get $x \ap z \ap u \ap x$ which implies $z \in [x]_{ap} \subset \mathcal{R}_{ap}$.

\end{pre}

The relation between being an attractor and being a quasi-attractor has been studied in the past, we recall this theorem of \cite{Ben98}.
\begin{prop}
Let $C$ be a non-empty subset of $\Delta$. The following assertions are equivalent :
\begin{enumerate}[(i)]
\item $C$ is an irreducible attractor i.e. it doesn't contain any proper attractor.
\item $C$ is an isolated quasi-attractor, i.e. there exists $U$, an open neighborhood of $C$, such that $U \cap \mathcal{R}_{ap} = C$.
\item $C$ is an isolated connected component of $\mathcal{R}_{ap}$ and \[C^+= \lbrace x \in M \, ; \, \exists y \in \mathcal{R}_{ap} \setminus C \text{ s.t. } C \ap y \ap x \rbrace = \emptyset\]
\end{enumerate}
\end{prop}

This result is proved in Part $5$ of \cite{Ben98}.

The following hypothesis is an analog of Hypothesis \ref{H1'} adapted to the context of ap-pseudo-orbits.

\begin{hyp}\label{h1}
There is only a finite number of ap-basic classes in $\mathring{\Delta}$ denoted by $K_i \, , \, i =1\cdots \eta$. We suppose that they are closed sets and indexed in such a way that the $k$ first $\lbrace K_i \rbrace_{i =1 \cdots k}$ are the quasi-attractors  and the $\eta-k$ others aren't.
\end{hyp}

\begin{prop}
For every $\theta$ small enough, there exist $\delta (\theta) \in [0 , \theta [$ and $T(\theta)\in ]0, \infty ]$ with $\delta( \theta), T(\theta)\not= (0 ,\infty)$ such that:

If there exists a $\delta , T$ pseudo-orbit $\xi_0 \cdots \xi_n$ verifying, $\delta < \delta (\theta)$ or $T > T(\theta)$ and, for a certain triplet $(i,i',j)\in \lbrace 1 \cdots \eta \rbrace^3$,
\[d(\xi_0 , K_i) < \delta \quad d(\xi_n , K_{i'}) < \delta  \quad d(\xi_j , K_i) > \theta\]
Then $i \not= i'$ and $K_i \prec K_{i'}$.
\end{prop}

This proposition means that, for $\delta$ small and $T$ large, the $\delta , T$ pseudo-orbits respect the partial order.

\begin{pre}

Suppose that, for every $\delta>0$ and every $T>0$, there exists a $\delta,T$-pseudo-orbit $\xi_0 \cdots \xi_n$ such that
\[d(\xi_0 , K_i) < \delta  \quad d(\xi_n , K_{i'}) < \delta \]
then, we can construct $\delta,T$-pseudo-orbits going from $K_i$ to $K_{i'}$ which in turn implies $K_i \ap K_{i'}$

Hence, if $K_i \not\ap K_{i'}$, there exists $\tilde{\delta} \geqslant 0$ and $\tilde{T} \in ]0 , + \infty ]$ such that a $\delta, T$ pseudo-orbit $\xi_0 \cdots \xi_n$ verifying 
\[d(\xi_0 , K_i) < \delta  \quad d(\xi_n , K_{i'}) < \delta \]
may only exist if $\delta > \tilde{\delta}$ or $T < \tilde{T}$.

Suppose now that $i=i'$, we will show that there exists $\hat{\delta}$ and $\hat{T}$ such that every $\delta, T$ pseudo-orbit  $\xi_0 \cdots \xi_n$ verifying either $\delta < \hat{\delta} $ or $T > \hat{T}$ and \[d(\xi_0 , K_i) < \delta  \quad d(\xi_n , K_{i}) < \delta \] doesn't contain any point at a distance greater than $\theta$ from $K_i$.

Suppose first that this assertion is false, thus we have real sequences $\delta_l \to 0$ and $T_l \to \infty$ and a sequence of $\delta_l,T_l$ pseudo-orbits verifying
\[d(\xi_0^l , K_i) < \delta \quad d(\xi_{n_l}^l , K_{i}) < \delta  \quad d(\xi_{j_l}^l , K_i) > \theta.\]

$\Delta$ being a compact set, we may suppose that, up to an extraction, as $l$ goes to infinity $x_0^l \to x \in K_i$, $x_{n_l}^l \to z \in K_i$ and $x_{j_l}^l \to y$ with $d(y,K_i) \geqslant \theta$. In that case we get $x \ap y \ap z$ and thus $y \in K_i$, which is absurd.

\end{pre}

\begin{prop}
For every $\delta >0$, there exists $T_0>0$ such that every $\delta, T_0$ pseudo-orbit intersects $N^\delta ( \mathcal{R}_{ap} )$.
\end{prop}

\begin{pre}

Let $x\in \Delta$ and $\gamma >0$, we define $T^{\gamma}(x) = Inf \lbrace t \geqslant 0 \,; \, \varphi_t (x) \in N^\gamma (\mathcal{R}_{ap}) \rbrace$

As $\omega(x) \cap \mathcal{R}_{ap} \not= \emptyset$ we get $T^\gamma (x) < + \infty$.

For $\alpha >0$ we will denote $N_\alpha= \lbrace x \in \Delta \, ; \, T^{\gamma} (x) \geqslant \alpha \rbrace$ the level sets of $T^\gamma$

Let us show that $N_\alpha$ is closed. Let $(x_n)_{n \in\N} $ be a sequence of elements of $N_\alpha$ converging to $y$. By the continuity of $\varphi_t$ we obtain then that, for every $t>0$, $\ds{\lim_{n \to \infty} \varphi_t(x_n)= \varphi_t(y)}$. If $t< \alpha$ we get $ \varphi_t(x_n) \in \left(N^\gamma (\mathcal{R}_{ap}) \right)^c$ which is closed. Hence we have $\varphi_t (y) \not\in N^\gamma (\mathcal{R}_{ap})$ for all $t<\alpha$ and thus $y \in N_\alpha$.

$T^\gamma (x)$ is then an upper semi-continuous function taking its values in $[0,+\infty [$,  $\Delta$ being a compact set, we know then that $T^\gamma(x)$ attains its maximum on $\Delta$ which we will denote $T^\gamma$.

Taking $T_0 > T^\delta$ gives us the result.

\end{pre}

The following corollary comes easily from the last two propositions.

\begin{cor}
For every $\delta >0$ and every $T>0$ there exists a family $V_i$ of open neighborhoods of the $K_i$ and positive real numbers $\delta_1$ and $T_1$ such that
\begin{itemize}
\item $\bar{N^{\delta_1} (K_i)} \subset V_i$ for $i=1 \cdots n$
\item Every $\delta_1, T_1$ pseudo-orbit starting from $V_i$ stays in $V_i$ for $i=1 \cdots k$
\item If there exists a $\delta_1 , T_1$ pseudo-orbit $x, (\xi_1,t_1) , \cdots , (\xi_p, t_p) , y$ with $x \in N^{\delta_1}(K_i)$ and $y \in N^{\delta_1}(K_j)$ such that 
\[\exists l \quad \exists\tilde{t} \in [0 ,t_l ] \text{ such that } \varphi_{\tilde{t}}(\xi_l) \not\in V_i\]
Then $i \not= i'$ and $K_i\prec K_{i'}$.
\item Every $\delta_1, T_1$ pseudo-orbit intersects $N^\delta ( \mathcal{R}_{ap} )$
\end{itemize}
\end{cor}

\begin{defin}
For $K$ compact subset of $\iD$ we denote
\[\beta_{\delta,K} (N) = \ds{ \sup_{x \in \Delta_N\cap K} \p_x[ \hat{X}^N (1) \in \Delta \setminus N^\delta (\varphi_1(x)) ]}\]
\end{defin}

\begin{prop}\label{bord}
If the flow $\lbrace \varphi_t \rbrace$ admits an attractor $A \subset \iD$, then, for all $K$ compact subset of $\iD$ and neighborhood of $A$, there exists $\delta >0$ such that $\rho_N^N \geqslant 1 - \beta_{\delta, K} (N)$.
Moreover, if there exists $V_K$ an open neighborhood of $\partial \Delta$ such that 
\[\ds{ \lim_{N \to \infty} \frac{\beta_{\delta,K} (N)}{\inf_{x \in V_K\cap \Delta_N} \p_x[\hat{X}^N (1) \in \partial \Delta ]}=0}\]
Then $\mu(V_K)=0$.
\end{prop}

\begin{pre}

Let $U$ be an open neighborhood of $A$ and $\delta >0$ such that $U \subset K$ and for every $t$, $\varphi_t(\bar{U}) \subset U$ and $N^\delta (\varphi_1(\bar{U})) \subset U$.

Then
\begin{align*} \rho_N^N \mu^N (U) & =  \sum_{x \in \Delta_N} p_N^N(x,U) \mu^N(x) \\
 & \geqslant  \sum_{x \in U\cap\Delta_N} \inf_{x \in U\cap\Delta_N} p_N^N(x,U) \mu^N(x) \\
 & \geqslant  \mu^N(U) \left( 1-\sup_{x \in U\cap\Delta_N} p_N^N(x,U^c)\right) \\
 & \geqslant  \mu^N(U) \left( 1-\sup_{x \in U\cap\Delta_N} p_N^N(x,(N^\delta (\varphi_1(\bar{U})))^c)\right) \\
 & \geqslant  \mu^N(U) \left( 1 - \beta_{\delta,K} (N) \right)\\
 \end{align*}
 
We finally get $\rho_N^N \geqslant 1-\beta_\delta (N)$

From this we obtain
\begin{align*} 1-\beta_{\delta,K} (N) & \leqslant  \rho_N^N \mu^N(\iD) \\
 & \leqslant  \sum_{x \in \Delta_N^*} \left( 1-p_N^N(x,\partial \Delta) \right) \mu^N(x) \\
 & \leqslant  \mu^N(M\setminus V_K) + \mu^N(V_K) \left( 1- \inf_{x \in V_0\cap\Delta_N} p_N^N(x,\partial \Delta) \right)\\
 \end{align*}
 
Hence
\[\mu^N(V_0) \leqslant \frac{\beta_{\delta ,K}(N)}{\inf_{x \in V_0\cap\Delta_N} p_N^N(x,\partial \Delta)}\]

$V_0$ being an open set, the weak convergence of the measures $\mu^N$ gives us the desired result.
\end{pre}

\begin{hyp}\label{h2}
We suppose that the flow $\varphi_t$ admits an attractor $A \subset \iD$ and that there exists $K$ a compact neighborhood of $A$ in $\iD$ and $V_K$ an open neighborhood of $\partial \Delta$ such that 
\[\ds{ \lim_{N \to \infty} \frac{\beta_{\delta,K} (N)}{\inf_{x \in V_K\cap\Delta_N} \p_x[\hat{X}^N (1) \in \partial \Delta ]}=0}\]
\end{hyp}

\vspace{3mm}

The following assumption is a technical one but we will see later that it is true under the first set of hypotheses.
\begin{hyp}\label{h3}
Let $j \in \lbrace k+1 , \cdots , \nu \rbrace$, i.e. such that $K_j$ is not a quasi-attractor.
Then  \[\exists \eta >0 \quad \forall \gamma >0 \quad \exists N_0 \quad \exists \chi : \N \to \R \quad \exists \zeta :\N \to \R \] such that  \[\lim_{n \to \infty} \zeta (n) = \lim_{n \to \infty} \chi (n)=0 \]
and
\[\sup_{x \in N^\eta (K_j)} \p_x \left[ \tau^N_{N^\eta (K_j)} > \e{\frac{\gamma}{\chi(N)}} \right] \leqslant \zeta(N)\]
\end{hyp}

We now arrive at the central theorem of this section.

\begin{thm}
Under Hypotheses \ref{h1}, \ref{h2} and \ref{h3}, the limiting measure $\mu$ has its support inside the union of ap-quasi-attractors $\ds{ \bigcup_{i=1}^k K_i}$
\end{thm}

\begin{pre}
Let $\delta>0$ small enough (how small will be specified later)


$K$ being an attractor, we know that $\rho_N \to 1$, that $\mu$ is $\varphi$-invariant and has its support inside $\mathcal{R}_{ap}$.

Hypothesis \ref{h2} and Proposition \ref{bord} gives us that $\mu(\partial \Delta)=0$, hence $\mu(K)=1$.

It only remains to be shown that, for every $j=k+1 \cdots \eta$, there exists an open neighborhood $W_j$ of $K_j$ such that $\mu(W_j)=0$

We know that
\[\rho_N^{N^2} \mu^N(W_j)  =  \sum_{x \in M} \p_x [ \hat{X}(N) \in W_j ]  = \int_M \p_x [ \hat{X}^N_N \in W_j ] d \mu^N(x)\]

We denote $t_N^i = \ds{ \left[ \frac{N}{i} \right]}$ and define the following events : 
\[E_N^\delta = \left\lbrace \exists \delta_1 , T_1 \text{ pseudo-orbit closer than }\delta \text{ from } \hat{X}^N_{t}, \; t \in [0,N] \right\rbrace \]
\[E'_N = \left\lbrace \forall i \in \lbrace k+1, \cdots , \eta \rbrace\: , \: \forall q \geqslant N^2 \: \hat{X}^N_p \in N^{\delta_1} (K_i) \Rightarrow \hat{X}^N_{p+q} \not\in N^{\delta_1}(K_i) \right\rbrace\]
$E'_N$ is the event "after its first entry in $N^{\delta_1}(K_i)$ the Markov chain $\hat{X}$ will have left it after $N^2$ steps".

Let $W_i$ be open neighborhoods of the $K_i$ such that $N^\delta (W_i) \subset V_i$ for every $i=1 \cdots \eta$. On $E_N^\delta \cap E'_N$ we get,

For $j \geqslant k$, 
\[X_{N}^N \in W_j \Rightarrow  \forall i \in \lbrace 2, \cdots b \rbrace \quad \hat{X}^N_{t_N^i}  \not\in N^\delta (K)\]

Hence 
\[\p [\hat{X}^N_{N} \in W_j \vert E_N^\delta \cap E'_N ] \leqslant \p \left[ \bigcap _{i=2}^b \lbrace  \hat{X}_{t^i_N}^N \not\in N^{\delta_1} (K) \rbrace\right]\]

On $E_N^\delta$, starting from $N^{\delta_1 - \delta}(K_i)$ the chain cannot enter $N^{\delta_1 - \delta}(K_i)$ once it has left $W_i^\delta$ (and also $V_i$).

Thus
\begin{align*}
\p [ (E'_N)^c \vert E_N^\delta ] & \leqslant  \p [ \exists i \text{ such that }\hat{X}\text{ doesn't leave }W_i\text{ before the time }N^2] \\
 & \leqslant  \sum_{i=k+1}^\eta \sup_{x \in N^{\delta_1-\delta}(K_i)} \p_x [ \hat{\tau}_{V_i}^N \geqslant N^2 ] \\
 & \leqslant  \sum_{i=k+1}^\eta \sup_{x \in N^{\delta_1-\delta}(K_i)} \p_x [ \hat{\tau}_{W_i}^N \geqslant N^2 ]
\end{align*}
where $\hat{\tau}_U^N = Inf \lbrace t \geqslant 0 \; ; \; \hat{X}_t^N \not\in U \rbrace$.

We now remark that, if $A,B$ and $C$ are three events, we get
\begin{align*}
\p[C] & =  \p[C \cap B^c \cap A ]+ \p[C \cap B \cap A ]+ \p[C \cap A^c ] \\
 & \leqslant  \p [ B^c \cap A ] + \p [ C \vert B \cap A ] + \p [A^c] \\
 & \leqslant  \p [B^c \vert A] + \p [ C \vert B \cap A ] + \p [A^c]
 \end{align*}
 
From that, we obtain
\[ \p_x[ \hat{X}^N_N \in W_j ] \leqslant \p_x[ (E_N^\delta)^c ] + \p_x [ (E'_N)^c \vert E_N^\delta ] + \p_x[\hat{X}^N_N \in W_j \vert E_N^\delta \cap E'_N ]\]

It only remains to control $\p_x[ (E_N^\delta)^c ]$.
In order to do that we will consider the pseudo-orbit $PO_N(t)$ defined by $x , (x,T_1) , (\hat{X}_{T_1}^N,T_1) \cdots $

For $N$ large enough, classic results on stochastic approximation algorithms (see e.g. \cite{BenWei03}) give us the following estimate :
\begin{align*}  \p_x[ (E_N^\delta)^c ] & \leqslant  \p [ PO_N(t) \text{ isn't a }\delta_1,T_1 \text{ pseudo-orbit} ] \\
& \qquad + \p_x[\sup_{t \in [0 , N]} d(\hat{X}^N_t,PO_N(t)) > \delta ]\\ 
& \leqslant  O((1-  \e{- \ep N})^{N})=O(\e{N \e{ -\ep N}})\end{align*}
with $\ep >0$

In conclusion we get
\begin{align*} \mu^N(V_j) & =  \frac{1}{\rho_N^{N}}\int_{\iD} \p_x[ \hat{X}^N_{N} \in W_j ] d\mu^N(x) \\
 &   \leqslant  \frac{1}{\rho_N^{N}}\int_{\iD} \p_x[ (E_N^\delta)^c ] + \p_x [ (E'_N)^c \vert E_N^\delta ] d\mu^N(x)\\
 & \qquad \quad +  \frac{1}{\rho_N^{N}}\int_{\iD}\p_x[\hat{X}^N_{N} \in W_j \vert E_N^\delta \cap E'_N ] d\mu^N(x) \\
 & \leqslant  \frac{1}{\rho_N^{N}} \left( O(\e{N \e{ -\ep N}}) + \zeta(N) + \sum_{i=1}^b \int_{\iD}  \p_x [ \hat{X}_{t_N^i} \not\in N^\delta (K) ] d\mu^N(x) \right) \\
 & \leqslant  \frac{1}{\rho_N^{N}} \left( O(\e{N \e{ -\ep N}}) + \zeta(N) + \sum_{i=1}^b   \rho_N^{N t_N^i} \mu^N( (N^\delta (K))^c )  \right) \\
 & \leqslant  \frac{1}{\rho_N^{N}} \left( O(\e{N \e{ -\ep N}}) + \zeta(N) + b \mu^N( (N^\delta (K))^c )  \right)
\end{align*} 
From Theorem \ref{attra}, we can infer that $\frac{1}{\rho_N^{N}} \to 1$.

Hence we have $\lim_{N \to \infty} \mu^N(V_j) =0$. As the sets $V_j$ are open neighborhoods of the sets $K_j$ we obtain $\mu(K_j)=0$.

\end{pre}

\subsection{Going back the hypotheses \ref{H1'}, \ref{H2'} and \ref{H3'}}
\hspace{1mm}\newline
The aim of this section is to prove that our first set of hypotheses implies the second one, thus proving the announced Theorem \ref{main}. In particular we will show that $\L$-quasi-attractors and ap-quasi attractors are the same.
\emph{In this section we will assume Hypotheses \ref{H1'}, \ref{H2'} and \ref{H3'} to be true.}
\begin{prop}
For every compact set $K\subset V_\alpha$, there exists a open neighborhood $V_K$ of $\partial \Delta$ such that
\[\lim_{N \to \infty} \frac{ \beta_{\delta,K}(N)}{\inf_{x \in V_K\cap \Delta_N} \p_x [\hat{X}^N(1) \in \partial \Delta ]}=0\]
where $\beta_{\delta,K}(N)=\ds{ \sup_{x \in K\cap \Delta_N} \p_x [ \hat{X}^N(1) \not\in N^\delta(\varphi_1(x))]}$
\end{prop}

\begin{pre}
Let $K$ be a compact set in $V_\alpha$ and let $\delta_0>0$ such that
\[ \delta_0 < \inf_{t \in [0,1],x \in K} d( \varphi_t(x), \partial \Delta) \]

Let $c(K)= \frac{1}{4} Inf \lbrace S(x,1,\phi) \vert x\in K , d(\phi,\varphi(x))\geqslant \delta_0 \rbrace >0$.

For $x \in K$ we have
\begin{align*} \limsup_{N \to \infty} \frac1N log \p_x  & \left[ \hat{X}^N \in \lbrace \phi \in \mathcal{C}_x([0,1]) ; d(\phi,\varphi(x)) >\delta_0 \rbrace \right] \\ & \leqslant  - \inf_{\phi ; d(\phi, \varphi(x)) > \delta_0} S(x,1,\phi)\\   & \leqslant  -4c(K)\end{align*}

Thus, for $N$ large enough and for every $x$ in $K$.
\[\frac1N log \p_x \left[ \hat{X}^N \in \left\lbrace \phi \in \mathcal{C}_x([0,1]) ; d(\phi,\varphi(x)) >\delta_0 \right\rbrace \right] \leqslant -3c(K)\]

Moreover, there exists $V_K$, an open neighborhood of $\partial \Delta$ such that
\[\lim_{N \to \infty} \inf_{x \in V_K} \frac1N log \p_x [\hat{X}^N(1) \in \partial \Delta ] \geqslant -2 c(K)\]

Then, for $N$ large enough, we have, for every $x$ in $V_K$
\[\p_x [\hat{X}^N(1) \in \partial \Delta ] \geqslant \e{- \frac{2c(K)}{1/N}}\]

Finally we get
\[\lim_{N \to \infty} \frac{ \beta_{\delta,K}(N)}{\inf_{x \in V_K} \p_x [\hat{X}^N(1) \in \partial \Delta ]} \leqslant \lim_{N \to \infty} \e{-\frac{c(K)}{1/N}} = 0\]
\end{pre}

\begin{prop}\label{contL}
Let $T>0$ and $K$ be a compact subset of $V_\alpha$.
Then, for every $\delta>0$,  there exists $\ep(K,T,\delta)$ such that
\[ \forall x \in K \quad S(x,T,\phi)\leqslant \ep \Rightarrow d(\phi,\varphi_.(x))_{[0,T]} \leqslant \delta\]
\end{prop}

\begin{pre}
Let's suppose that this result is false, then
\[\exists \delta \quad \forall \ep \quad \exists x \quad \exists \phi_\ep \text{ s.t. } S(x,T,\phi_\ep) < \ep \text{ and }d(\varphi_.(x),\phi_\ep) > \delta\]

We also know that $c= \inf \lbrace S(x,T, \phi) \; ; x \in K \; ; \; d(\phi,\varphi_.(x))> \delta \rbrace >0$

Thus, for every $\ep >0$ 
\[\ep > S(X,T, \phi_\ep ) \geqslant c >0\] which is absurd. 
\end{pre}

\begin{prop}
The function $B_\L : V_\alpha \times V_\alpha \to \R_+$ is upper semi-continuous.
\end{prop}

\begin{pre}
Let $x_n$ and $y_n$ two sequences in $V_\alpha$ converging to $x$ and $y$ in $V_\alpha$

For every $\delta >0$ there exists a path $x=\xi_0, \cdots \xi_{n(\xi)}=y$ such that
\[B_\L (x,y) \leqslant A_{n(\xi)}(\xi) \leqslant B_\L(x,y)+\delta \]

Let's now consider the path $\xi^n$ given by $x_n=\xi_0^N , \xi_1 ,\cdots \xi_{n(\xi)-1} , y_n=\xi_{n(\xi)}$

If we take $n$ such that 
\[\vert \L(x_n,\xi_1)-\L(x,\xi_1) \vert < \delta \text{ and }\vert \L(\xi_{n(\xi)-1},y_n)-\L(\xi_{n(\xi)-1},y) \vert < \delta\]
We obtain
\[B_\L(x_n,y_n) \leqslant A_{n(\xi)}(\xi^N) \leqslant B_\L(x,y) +3\delta\]

Thus, for every $\delta>0$, $\ds{\limsup_{n \to \infty} B_\L(x_n,y_n)}    \leqslant B_\L(x,y) +3\delta$

Finally  $\ds{\limsup_{n \to \infty} B_\L(x_n,y_n)}  \leqslant B_\L(x,y)$
\end{pre}

\begin{prop}\label{proche}
Let $x \in \mathcal{R}_\L$ such that $[x]_\L$ is a closed subset of $V_\alpha$. Let $\theta>0$ such that $N^\theta([x]_\L) \subset V_\alpha$.

Then,\[ \exists \delta>0 \quad \exists T>0 \quad \forall \Psi \text{ with } \Psi(0) \in  [x]_\L, \quad \Psi(\tilde{T}) \in [x]_\L, \quad \tilde{T}>T,\]
\[ \text{ and } S(\Psi(0),\tilde{T},\Psi)<\delta \quad \forall t\in [0,\tilde{T}] \qquad \Psi(t) \in N^\theta([x]_\L).\]
\end{prop}

\begin{pre}
Let $x \in \mathcal{R}_\L$ such that $[x]_\L$ is a closed subset of $V_\alpha$ and let  $\theta>0$ such that $N^\theta([x]_\L) \subset V_\alpha$.

Suppose that the announced result is false, then, there exists a family of functions $\Psi_n\in \mathcal{C}([0,T_n],V_\alpha)$ such that  
\[\Psi_n(0) \in  [x]_\L \quad \Psi_n(T_n) \in [x]_\L \quad \lim_{n \to \infty} T_n =\infty \quad S(\Psi(0),T_n,\Psi)<1/n \] and
\[ \forall n \quad \exists t_n>0 \text{ such that }\Psi_n(t_n) \not\in N^\theta([x]_\L).\] 

$[x]_\L$ being a compact set, we can assume without loss of generality that $\Psi_n(0) \to u \in [x]_\L$ and $\Psi_n(T_n) \to v \in [x]_\L$.

Let $\tau_n=\inf \lbrace t>0 \,;\, \Psi_n(t) \not\in N^\theta([x]_\L )\rbrace$. Then we have \[\Psi_n(\tau_n) \in \adh{V_\alpha} \setminus N^\theta([x]_\L)\] 
and $ \adh{V_\alpha} \setminus N^\theta([x]_\L)$ is a compact set.

Thus, without loss of generality, we can assume that \[\Psi_n(\tau_n) \to w \in V_\alpha \setminus N^\theta([x]_\L).\]

Using the sequential continuity of $\L$  and the fact that $u$ and $v$ are $\L$-chain-recurrent points, we get that $u \leadsto_\L w$ and $w \leadsto_\L v$. Thus $w\in [x]_\L$ and we get a contradiction

\end{pre}

\begin{cor}\label{prochebis}
Suppose $\mathcal{R}_\L$ is a closed subset of $V_\alpha$. Let $\theta>0$ such that $N^\theta(\mathcal{R}_\L) \subset V_\alpha$
Then, 
\[ \exists \delta>0 \quad \exists T>0 \quad \forall \Psi \text{ with }\Psi(0) \in  \mathcal{R}_\L ,\quad \Psi(\tilde{T}) \in \mathcal{R}_\L, \quad \tilde{T}>T \]
\[\text{ and }S(\Psi(0),\tilde{T},\Psi)<\delta \quad \forall t\in [0,\tilde{T}] \qquad \Psi(t) \in N^\theta(\mathcal{R}_\L). \]
\end{cor}

\begin{prop}\label{pas1}
Let $x \in \mathcal{R}_\L$ such that $[x]_\L$ is a closed subset of $V_\alpha$. Then $x \in \mathcal{R}_{ap}$ and $[x]_\L \subset [x]_{ap}$.
\end{prop}

\begin{pre}
Let $y \in [x]_\L$.
Let $(T_n)_{n \in \N}$ be a sequence of positive real numbers and $\Psi_n$ be a sequence of functions from $[0,T_n]$ to $V_\alpha$ such that $\Psi_n(0) = x$, $\Psi_n(T_n) =y$ and $S(x,T_n,\Psi_n)<1/n$.

$T_n$ is either a bounded sequence or it goes to infinity (up to a sub-sequence).

Let's suppose that $(T_n)_{n \in \N}$ is an increasing sequence of positive real numbers going to infinity and $\Psi_n$ a sequence of functions such that $\Psi_n(0) = x$, $\Psi_n(T_n) =y$ and $S(x,T_n,\Psi_n)<1/n$

The former proposition gives us a compact set $K\subset V_\alpha$, and positive numbers $\delta$ and $T$ such that, for $T_n>T$ and $1/n <\delta$, $\Psi_n$ lives in $K$

Let $\tilde{T}>T$, $\tilde{\delta}<\delta$ and $\ep=\ep(\tilde{\delta},\tilde{T},K)$ given by Proposition \ref{contL}.

Then, from Proposition \ref{contL}, we can infer that, for $n$ large enough, \[d(\Psi_n(\tilde{T}),\varphi_{\tilde{T}}(x)) \leqslant \delta\]

As $\Psi_n(\tilde{T}) \in K$, we also have \[d(\varphi_{\tilde{T}}(\Psi_n(\tilde{T})),\Psi_n(2\tilde{T}))< \delta\]

By iterating this process we get that, for $n$ large enough, \[(x,\tilde{T}),(\Psi_n(\tilde{T}),\tilde{T}), \cdots \] is a $\delta, \tilde{T}$ ap-pseudo orbit linking $x$ to $y$.

Hence $x \ap y$.

\vspace{3mm}

Let's suppose now that $T_n$ is bounded by $\bar{T}$. By taking a sub-sequence we can assume that $T_n \to T$ as $n \to \infty$ for some $T$ with $\bar{T} \geqslant T \geqslant 0$.  we continue the function $\Psi_n$ to the time $\bar{T}$ by concatenating the flow $\varphi_{\centerdot}$. Then we know that $\Psi_n \to \varphi_\centerdot (x)$ uniformly over $[0,\bar{T}]$ as $n \to \infty$. Hence $y=\lim_{n \to \infty} \Psi_n(T_n)=\varphi_T(x)$. Thus $T>0$ and $y \in \gamma^+(x)$. 

If $y=x$ we get then that, either $x \ap x$ or $x$ is periodic which implies $x \in \mathcal{R}_{ap}$.

Then, if $y \in [x]_\L$ we obtain that either $x \ap y$ or $y \in \gamma^+(x) \subset [x]_{ap}$. Taking $y=x$ gives us $x\in \mathcal{R}_{ap}$, i.e. $\mathcal{R}_{\L} \subset \mathcal{R}_{ap}$. As the roles of $x$ and $y$ can be exchanged we also get $y \ap x$ and in conclusion $[x]_\L \subset [x]_{ap}$.

\end{pre}

\begin{prop}
Closed $\L$-classes in $V_\alpha$ are positively invariant sets for the flow $\varphi$
\end{prop}

\begin{pre}
Let $x \in \mathcal{R}_\L$ such that $[x]_\L$ is a closed $\L$ basic class in $V_\alpha$. 
Let $T>0$ and $\theta$ such that $N^\theta([x]_\L) \subset V_\alpha$.

If the paths linking $x$ to itself have bounded length then $x$ is periodic and in this case $\gamma^+(x) \subset [x]_\L$.

We now suppose that the paths have unbounded length.

Let $\ep=\ep([x]_\L,\tau,\theta)$ and let $\delta>0$ and $\tau>0$ given by Proposition \ref{proche}, let $\Psi$ be a path of length greater than $\tau$ and of cost smaller than $min(\ep,\delta)$ linking $x$ to itself. 

Then, by the triangular inequality we get
\[d(\varphi_T(x),[x]_\L) \leqslant d(\varphi_T(x),\Psi(T)) + d(\Psi(T), [x]_\L) \leqslant 2 \theta\]

Thus, by making $\theta$ go to zero, we obtain $\varphi_T(x) \in \bar{[x]_\L} =[x]_\L$;
\end{pre}

\begin{lem}
Let $\rightarrowtail$ be a binary, transitive relation on $V_\alpha$, such that, for every $x \in V_\alpha$, every $y\in \alpha(x)$ and every $z \in \omega(x)$,
\[y\rightarrowtail x \quad x \rightarrowtail z \quad z\rightarrowtail z\]
Let $[x]_\rightarrowtail$ be a maximal, closed and isolated $\rightarrowtail$-basic class.

Then $[x]_\rightarrowtail$ is an attractor for the flow $\varphi$
\end{lem}

\begin{pre}
Let $W$ be an isolating open neighborhood of $[x]_\rightarrowtail$, i.e. $W$ is an open neighborhood of $[x]_\rightarrowtail$ such that the only $\rightarrowtail$-recurrent points of $W$ belong to $[x]_\rightarrowtail$.

If $[x]_\rightarrowtail$ isn't an attractor then there exists $p \in \partial \bar{W}$ such that
\[\gamma_-(p) \subset W \text{ and } \alpha(p) \subset [x]_\rightarrowtail\]

Let $y \in \omega(p)$, we have $p\rightarrowtail y$ and $y\rightarrowtail y$. 

Similarly let $z\in \alpha(p)$, we have $z\rightarrowtail p$ and thus $z \rightarrowtail y$.

Let's define 
\[[x]_\rightarrowtail^+ = \lbrace z \; ,  \exists y \in \mathcal{R}_\rightarrowtail \setminus [x]_\rightarrowtail  \text{ and } x \rightarrowtail y \rightarrowtail z \rbrace\]

As $[x]_\rightarrowtail$ is maximal, we have $[x]_\rightarrowtail^+ =\emptyset$.

Yet we have $x\rightarrowtail y$ and $y \in \mathcal{R}_\rightarrowtail$.

Hence $y \in [x]_\rightarrowtail$, and, as $x \rightarrowtail p\rightarrowtail y$ we get $p \in [x]_\rightarrowtail$ which is absurd.
\end{pre}

\begin{prop}
Let $[x]_\L$ be a closed isolated quasi-attractor. Then $[x]_\L$ is an attractor.
\end{prop}

\begin{prop}\label{2.15}
Let $K$ be a compact subset of $V_\alpha$.

For $\eta >0$, $\delta>0$, $\tilde{T} >0$ there exists $N_0 >0$ such that, for all $N >N_0$, all $x\in K$, all $T < \tilde{T}$ and all $\Psi \in \mathcal{C}_x([0,T],K)$
 \[\p_x[ d(\Psi , \hat{X}_N ) < \eta ] \geqslant exp \left( -\frac{S(x,T, \Psi) + \delta}{1/N} \right)\]
\end{prop}

\begin{pre}
Let $\xi_{\gamma,T}^K = sup \lbrace \vert S(x,T,\Psi) - S(x,T,\phi) \vert \, ;\, x,y \in K \, , \, d(\Psi,\phi)< \gamma \rbrace$.
The inferior semi-continuity of $S$ gives us that $\ds{ \lim_{\gamma \to 0} \xi_{\gamma,T}^K =0 }$

Let $\eta >0$, $\delta>0$, $\tilde{T} >0$ and let $\gamma < \eta$ such that $\xi_{\gamma, \tilde{T} }^K < \delta $ and $N^\gamma (K) \subset V_\alpha$

The large deviation principle gives us that, for $\phi \in \mathcal{C}_x([0,T],K)$,
\[\liminf_{N \to \infty} \frac1N log \inf_{x\in K}\p_x [ \hat{X} \in N_x^\gamma (\phi) ]]\geqslant - \sup_{x\in K} \inf_{\Psi \in N^\gamma(\phi)} S(x,T, \Psi) \]
where $N^\gamma(\phi)=\lbrace \Psi \in \mathcal{C}([0,T],K) \; ; \; \Vert \Psi - \phi \Vert < \gamma \rbrace$.
Thus, there exists $g_K(N)$ a function which goes to zero as $N$ goes to infinity such that

\[  \frac1N log  \inf_{x\in K}\p_x [ \hat{X} \in N_x^\gamma (\phi) ]]\geqslant -  \sup_{x\in K}\inf_{\Psi \in N^\gamma(\phi)} S(x,T \Psi)  -g_K(N)\]

Hence 
\begin{align*} \quad\p_x[ d(\Psi , \hat{X}_N ) < \eta ]  &\geqslant  \quad\p_x[ d(\Psi , \hat{X}_N ) < \gamma ] \geqslant  \inf_{x\in K} \p_x [ \hat{X} \in N^\gamma (\Psi)]\\ & \geqslant  exp \left( \frac{ - \sup_{x\in K}\inf_{\phi \in N^\gamma(\Psi)} S(x,T,\phi) -g_K(N)}{1/N}\right)\\ & \geqslant  exp \left( \frac{ -\xi_{\gamma,T}^K  - S(x,T,\Psi) -g_K(N)}{1/N}\right)
\end{align*}

Taking $N$ large enough such that $\xi_{\gamma,T}^K +g_K(N) <\delta$ gives us the result
\end{pre}

\begin{prop}
We denote by $BC(\varphi)=\bar{ \lbrace x \in \Delta \, ; \, x \in \omega (x) \rbrace}$ the Birkhoff center of the flow $\varphi$. Then \[V_\alpha \cap BC(\varphi) \subset \mathcal{R}_\L\]
\end{prop}

\begin{pre}
We know that $BC(\varphi)=\bar{ Rec(\varphi)}=\bar{ \lbrace x \in \Delta \, ; \, x \in \omega (x) \rbrace}$. It's apparent that $Rec(\varphi)\subset \mathcal{R}_\L$. Hypothesis \ref{H3'} allow us to conclude.
\end{pre}

\begin{cor}
Let $\mu$ be an invariant measure for the flow $\varphi$ whose support $S$ lies within $V_\alpha$. Then $S \subset \mathcal{R}_\L$.
\end{cor}

We arrive at the main theorem of this section, linking $\L$-chain recurrence with ap-chain recurrence.

\begin{thm}\label{main2}
Under the hypotheses \ref{H1'}, \ref{H2'} and \ref{H3'} we have 
\[\mathcal{R}_{ap} \cap V_\alpha = \mathcal{R}_\L \cap V_\alpha\]  
\[ \forall x \in \mathcal{R}_\L \cap V_\alpha \qquad [x]_\L=[x]_{ap} \]
\end{thm}

\begin{pre}
By Proposition \ref{pas1} we already know that, if $[x]_\L$ is a closed $\L$-basic class in $\iD$, then $x\in \mathcal{R}_{ap}$ and $[x]_\L \subset [x]_{ap}$.

The function $B_\L$ is upper semi-continuous. Thus
\[\forall \gamma >0 \quad \exists U_\gamma^i \text{ neighborhood of }K_i \text{ s.t. } \forall (a,b)\in (U_\gamma^i)^2 \quad B_\L(a,b)<\gamma\]

\begin{lem}\label{lem1}
Let $U$ be a neighborhood of $\mathcal{R}_\L \cap \iD$ such that $K=\bar{U}$ is a compact subset of $\iD$. Thus there exists $T>0$ such that, if $\gamma^+(x) \subset K$ then $\lbrace \varphi_t(x) \, ; \, t \in [0,T] \rbrace \cap U \neq \emptyset$  
\end{lem}

\begin{pre}
Due to the fact that, if $\gamma_+(x) \subset K$, $\omega(x) \subset U$ we can't have $\gamma^+(x) \subset K \setminus U$.
Let's define $\nu(x) = inf \lbrace t \, ; \, \varphi_t(x) \not\in K \setminus U \rbrace < +\infty$

We will show that $\nu(x)$ is upper semi-continuous.

Let $x \in K$ and $x_n \in K^\N$ such that $x_n \to x$. 

The continuity of $\varphi$ gives us $d(\varphi_{T_n}(x_n),\varphi_T(x)) \to 0$ with the convention that $\varphi_{\infty}(x)=\omega(x)$. Let $\tilde{T}=\nu(x)$, there exists a sequence $\ep_m \to 0$ such that $\varphi_{\tilde{T} + \ep_m}(x) \not \in K \setminus U$.

Yet we have $\varphi_{\tilde{T} + \ep_m}(x_n) \to \varphi_{\tilde{T} + \ep_m}(x)$.
Hence, for $n$ large enough we get $\nu(x_n) \leqslant \tilde{T}+\ep_m$

We now get $\ds{\limsup_{n \to \infty} \nu(x_n) \leqslant \tilde{T} + \ep_m}$

Then $\ds{\limsup_{n \to \infty} \nu(x_n) \leqslant \nu(x)}$, i.e. $\nu$ is u.s.c.

Thus $T=\ds{\max_{x \in K\setminus U} \nu(x)}$ exists and we get, if $\gamma^+(x) \subset K$ then \[\lbrace \varphi_t(x) \, ; \, t \in [0,T+\ep] \rbrace \cap U \neq \emptyset\]

\end{pre}
%
%
%
%

Let $x \in \mathcal{R}_{ap} \cap \iD$ and let $y \in [x]_{ap}$, there exists two sequences of positive real numbers $\delta_k \to 0$, $T_k \to \infty$ and a sequence of $\delta_k , T_k$ ap-pseudo-orbits linking $x$ to $y$ denoted $\Psi^k$.

Lemma \ref{lem1} gives us $K\subset \iD$ a compact neighborhood of $\mathcal{R}_{ap}$ such that, for $k$ large enough, every $\delta_k ,T_k$ ap-pseudo-orbit stays in $K$.
 
We know that $\mathcal{R}_\L \subset \mathcal{R}_{ap}\subset K$. Let $\gamma$ such that $\ds{\bigcup_{i} N^\gamma (K_i) \subset K}$. Let then $T$ given by the former lemma. For $k$ large enough, we have $T_k>T$ and thus every continuous part of the $\delta_k ,T_k$ ap-pseudo-orbit intersect $U_\gamma$.

Let $\ep < \gamma$ and $V= \ds{\bigcup_i N^{\ep}(K_i)}$, we assume $\ep$ to be small enough such that the $\ep$-neighborhood of the $K_i$ are disjoint.

We define two sequences of times $(\sigma_i(k))_{i=1 \cdots p_k}$ and $(\tau_i(k))_{i=1 \cdots q_k}$ by
\[ \sigma_0(k)=0 \qquad \tau_0(k)=min \lbrace t \, ; \, \Psi^k_t \not\in V \rbrace\]
\[ \sigma_{i+1}(k) = min \lbrace t > \tau_i(k) \, ; \, \Psi^k_t \in V \rbrace \quad \tau_{i+1}(k) = min \lbrace t > \sigma_{i+1}(k) \, ; \, \Psi^k_t \not\in V \rbrace\]

By Lemma \ref{lem1} we know that $\sigma_{i+1}(k) -\tau_i(k) < 2 T $.

If our pseudo-orbit $\psi^k$ enters more than once the same $N^\ep( K_i)$ then we truncate what happens between the first entry and the last exit and we will keep the same name for the new path (which may no more be a pseudo-orbit). This way we get $q_k , p_k < \eta$. 

If $y \in \ds{\bigcup_i K_i}$ then $q_k = p_k+1$, else $q_k=p_k$.

In the first case we get
\[B_\L(x,y) \leqslant \sum_{i=0}^{q_k-1} B_\L ( \Psi_{\sigma_i(k)}^k, \Psi_{\tau_i(k)}^k) + \sum_{i=0}^{p_k} B_\L ( \Psi_{\tau_i(k)}^k, \Psi_{\sigma_{i+1}(k)}^k) +B_\L ( \Psi_{\tau_{q_k}}^k ,y)\]

In the second case, we get 
\[B_\L(x,y) \leqslant \sum_{i=0}^{q_k-1} B_\L ( \Psi_{\sigma_i(k)}^k, \Psi_{\tau_i(k)}^k) + \sum_{i=0}^{p_k} B_\L ( \Psi_{\tau_i(k)}^k, \Psi_{\sigma_{i+1}(k)}^k) +B_\L ( \Psi_{\sigma_{p_k}}^k ,y)\]

In either case we obtain
\[B_\L (x,y) \leqslant (1+ \nu) \alpha + (1+\nu) Sup \left\lbrace B_\L (a,b) \; ; \; d(b , \varphi_{[0,2T]}(a)) < \ep_k \right\rbrace\]
where \[\ep_k=\sup_{d(x,y)< \delta_k\, , \, t\in [0,2T]} d(\varphi_t(x),\varphi_t(y))\]

From the continuity of $\varphi$, the positive invariance of $[x]_{ap}$ and Proposition \ref{cont} we infer that
\[ \lim_{k \to 0} Sup \left\lbrace B_\L (a,b) \; ; \; d(b , \varphi_{[0,2T]}(a)) < \ep_k \right\rbrace =0\]

Hence we have, for every $\alpha>0$, $B_\L(x,y) \leqslant (\nu +1) \alpha$.

In conclusion $x \leadsto_\L y$. Similarly $y \leadsto_\L x$, i.e. $y \in [x]_\L$

\end{pre}

\textit{\underline{Remark :}} 
This proof gives us a little more, it proves that, if $y\in \mathcal{R}_{ap}$ and $x \ap y$, then $x\leadsto_\L y$.

\vspace{3mm}

When we take $\alpha$ such that $\iD \cap \mathcal{R}_\L = V_\alpha \cap \mathcal{R}_\L$, this proposition proves that the first set of hypotheses implies Hypothesis \ref{h3}.
\begin{prop}

Let $j \in \lbrace k+1 , \cdots , \eta \rbrace$, i.e. such that $K_j$ is not a ap-quasi-attractor.
Then  \[\exists \lambda >0 \quad \forall \gamma >0 \quad \exists N_0  \quad \exists \zeta :\N \to \R \text{ with } \lim_{n \to \infty} \zeta (n) =0 \text{ such that}\]
\[\sup_{x \in N^\lambda (K_j)} \p_x \left[ \tau^N_{N^\lambda (K_j)} > \e{\frac{\gamma}{1/N}} \right] \leqslant \zeta(N)\]
\end{prop}

\begin{pre}
$K_j$ isn't a quasi-attractor, thus $\exists \lambda>0$ such that $\bar{N^{2 \lambda}(K_j)} \subset \iD$ and, for all $\gamma >0$ and all $x \in N^\lambda(K_j)$, there exists $T^\gamma$ and $\Psi^\gamma$ such that
\[ \Psi^\gamma(0)=x, \quad \Psi^\gamma(T^\gamma)=y^\gamma \not\in N^{2\lambda}(K_j), \quad S(x, T^\gamma,\Psi^\gamma) < \gamma \]

Let then $U=N^{2\lambda}(K_j)$, $V=\omega(\bar{U})$. Let $r>0$ and let $K$ be a compact subset of $\iD$ such that $N^r(U)\subset K$.

As $\Psi^\gamma$ starts in $U$ and ends outside of $U$ it must pass through $K \setminus U$. Without loss of generality we may suppose that $\Psi^\gamma$ stays in $K$ and $y^\gamma \in K \setminus U$.

Proposition \ref{2.15} says that :

For $\lambda >0$, $\delta>0$, $\tilde{T} >0$ , there exists $N_0(K) >0$ such that, for all $N >N_0$, all $T < \tilde{T}$, all $x\in K$ and all $\Psi \in \mathcal{C}_x([0,T],K)$ \[\p_x[ d(\Psi , \hat{X}_N ) < \lambda ] \geqslant exp \left( -\frac{S(x,T, \Psi) + \delta}{1/N} \right)\]

Applying this to our case gives us $N_0$ such that, $\forall N >N_0$
\[ \p_x[d( \Psi^\gamma, \hat{X}^N ) < \lambda ] \geqslant exp\left( -\frac{S(x,T^\gamma, \Psi^\gamma) +\delta}{1/N} \right) \geqslant exp \left( - \frac{\gamma +\delta}{1/N} \right)\]

Yet $\lbrace d(\hat{X}^N,\Psi^\gamma) < \lambda \rbrace$ implies that $\hat{X}^N$ leaves $N^\lambda(K_j)$ before the time $T^\gamma$

Hence \[\forall x \in K \quad \forall N \geqslant N_0(K) \quad\p_x [ \tau_{N^\lambda(K_j)} > T^\gamma ] < 1-\e{-\frac{\gamma + \delta }{1/N}}\]

Then
\begin{align*}
\p_x \left[ \tau_{N^\eta(K_j)} > \e{\frac{\gamma}{1/N}} \right] & <  \left( 1-\e{-\frac{\gamma + \delta }{1/N}} \right)^{\left[ \frac{\e{\frac{2\gamma}{1/N}}}{T^\gamma} \right]} \cr
 & <  \e{\left[ \frac{\e{\frac{2\gamma}{1/N}}}{T^\gamma} \right] ln\left( 1-\e{-\frac{\gamma + \delta}{1/N}} \right)} \cr
 & <  \e{\left[ \frac{\e{\frac{2\gamma}{1/N}}}{T^\gamma} \right] -\e{-\frac{\gamma + \delta}{1/N}} } \cr
 & <  \e{-\e{-\frac{\gamma -\delta}{1/N}}\frac{1}{T^\gamma}} \end{align*}
 
Taking $\delta < \gamma$ and $\zeta(N) =   \e{-\e{-\frac{\gamma -\delta}{1/N}}}$ allows us to conclude.

\end{pre}

We linked ap-basic classes with $\L$-basic classes but, unless the partial orders on the class are similar, we might not have the same quasi-attractors for both partial orders. This proposition shows that the quasi-attractors are the same.

\begin{prop}
Let $[x]$ be a basic-class. $[x]$ is a $\L$-quasi-attractor if and only if $[x]$ is an ap-quasi-attractor.
\end{prop}

\begin{pre}

Theorem \ref{main2} already gives us that, if $y\in \mathcal{R}_{ap}$ and $x \ap y$, then $x\leadsto_\L y$.

Hence, if $[x]$ is a $\L$-quasi-attractor, it is also an ap-quasi-attractor.

Suppose now that $[x]$ is an ap-quasi-attractor.

Let $y \in \mathcal{R}_\L$ such that $x \leadsto_\L y$. If the path linking $x$ to $y$ have bounded length we get $y \in \gamma^+(x)$. As $\L$-basic classes are positive invariant sets we get then $y\in[x]$.

Let us now suppose that the paths linking $x$ to $y$ have unbounded length. Propositions \ref{prochebis} and \ref{cont} give us a compact set containing the paths linking $x$ to $y$ for $T$ large enough. Using the same technique as in the proof of the proposition \ref{pas1} gives us $x \ap y$, i.e. $y\in [x]$.
\end{pre}

\appendix
\section{Reminders}
\hspace{1mm}\newline
The purpose of this appendix is to give small reminders about some basic properties of the objects used here. We won't give extensive proofs of the results or complete theory of these objects but will include some references should the reader want to know more.

\subsection{Quasi-stationary distributions}
\hspace{1mm}\newline
This subsection has been inspired by the notes of a course given by Sylvie M\'el\'eard at the VI Escuela de Probabilidad y Procesos Estoc\`asticos in Guanajuato in September 2009.

Here we will consider a Markov chain $Z_t$, either in discrete or continuous time, whose state space $E \subset \R^d$ admits an absorbing state, denoted by $\lbrace 0 \rbrace$. We will denote $E^*=E \setminus \lbrace 0 \rbrace$, $\mathcal{P}^*$ the set of probability measures whose support lies in $E^*$. We define $T_0$ to be the absorption time. 
\[T_0=\inf \lbrace t>0 ; Z_t=0 \rbrace\]
We suppose that, whatever the initial state is, the process will almost surely be absorbed, i.e.
\[ \forall z \in E \qquad \p_z [T_0 < \infty ]=1\]

\begin{defin}
\hspace{0.1mm} \newline
\begin{enumerate}
\item A probability measure $\mu$ on $E^*$ is said to be a quasi-stationary distribution (QSD) if and only if, for every Borelian set $A \subset E^*$ and every $t>0$,
\[ \p_\mu [Z_t \in A \vert T_0 > t ] = \mu(A)\]
We remark that, in this case, $\mu$ is a fixed point for the conditional evolution 
\[\nu \mapsto \p_\nu [ Z_t \in \centerdot \vert T_0 > t ]\]
\item A quasi-limiting measure (QLD) for $\alpha \in \mathcal{P}^*$ is a probability measure $\nu$ on $E^*$ such that, is the weak convergence sense,
\[  \quad \lim_{t \to \infty} \p_\alpha [ Z_t \in \centerdot \vert T_0 >t ]=\nu\]
\item The Yaglom limit is the probability measure $\pi$ on $E^*$ defined by, for $A \in \mathcal{B}(E^*)$,
\[\pi(A) = \lim_{t \to \infty} \p_z [ Z_t \in A \vert T_0 >t ]\]
as soon as this limit exists and doesn't depend on $z \in E^*$
\end{enumerate}
\end{defin}

\begin{prop}
Suppose that $\mu$ is a QSD for this process $Z_t$. Then there exists a positive real number $\theta(\mu)$ such that 
\[\p_\mu [T_0 > t ] = \e{- \theta(\mu)t}\]
\end{prop}

See \cite{MelVil11} for a proof

The following theorem is true in a broader context (see e.g. Sylvie M\'el\'eard course notes or Bonsall paper \cite{Bon57}) but this version is sufficient for the problem at hand.

\begin{thm}[Perron-Frobenius Theorem]\label{QSD}

Let $E$ be a finite set and let $X_n$ be a discrete time Markov chain on $E$ with transition matrix $Q$. Up to merging some states, we suppose that $X_n$ admits an unique absorbing state $\lbrace 0 \rbrace$. We denote by $Q^*$ the matrix $Q$ restricted to $E^*=E \setminus \lbrace 0 \rbrace$. We suppose $Q^*$ to be irreducible. Then
\begin{enumerate}[1)]
\item There exists an unique quasi-stationary distribution (QSD) obtained as the only left eigenvector $\mu$ of $Q^*$ verifying
\[\mu_i >0 \quad ; \quad \sum_i \mu_i =1\]
The corresponding eigenvalue $0 <\rho= \e{-\theta}<1$ is such that
\[\p_{\mu} [T_0 > n ] = \e{- \theta n} \] 
\item The measure $\mu$ is the Yaglom limit, i.e. 
\[ \forall i,j \in E^* \lim_{n \to \infty} \p_i [X_n =j \vert T_0 >n ]=\mu_j\]
\item For every couple $(i,j)$ of points in $E^*$ $\displaystyle{ \lim_{n \to \infty} \e{\theta n} \p_i [X_n=j]=\pi_i \mu_j}$
where $\pi$ is the unique right eigenvector of $Q^*$ associated to $\rho$ and such that $\pi_i >0 \quad ; \quad \sum_i \mu_i \pi_i =1$
\item For every couple $(i,j)$ of $E^*$ and every $n,m \in \N^*$  \[\displaystyle{ \lim_{n \to \infty} \frac{\p_i[T_0 >n+m]}{ \p_i [T_0>n]}=\frac{\pi_i}{ \pi_j} \e{- \theta n}}\]
\end{enumerate}
\end{thm}

The result concerning the Yaglom limit can be refined

\begin{prop}
We suppose $Q$ to be $\C$ diagonalizable (this will be the case most of the time as the set of diagonalizable matrix of size $n$ contains an open subset that is dense in $\mathcal{M}_n(\C)$). We denote  $\rho= \lambda_1 > \vert \lambda_2 \vert > ... \geqslant \vert \lambda_n \vert$ the eigenvalues of $Q$. Let $\nu$ be a measure on $E$, then :
\[\Vert \p_\nu [X_k \in \centerdot \vert T_0 >k ] - \mu \Vert_2 \leqslant C \left(\frac{\vert \lambda_2 \vert}{\rho}\right)^k\]
\end{prop}

This result is an interpretation in our context of a classic result in numerical analysis on the power method, see e.g. \cite{MR1751750}.

For more information on the subject see \cite{MelVil11} and \cite{Pol12}.

\subsection{Stochastic approximation algorithms}
\hspace{1mm}\newline
This subsection has been inspired by \cite{BenWei03} where you can find more information on the subject and complete proofs of the results.

Here we will consider a family of discrete processes $(X_n^N)_{n \in \N}$ taking values in $K$ a compact subset of $\R^m$. The parameter $N$ indexing the processes may take either real or integer values. We suppose that the $X_n^N$ are defined on a probability space $(\Omega, \mathcal{F}, \mathbb{P})$ and denote by $\mathcal{F}_n^N$ the $\sigma$-algebra generated by $\lbrace X_i^N, i=1,...,n \rbrace$. We suppose that
\[ X_{n+1}^N -X_n^N = \frac{1}{N} (F(X_n^N)+U_{n+1}^N)\]
where
\begin{enumerate}[(i)]
\item $F: \R^m \to \R^m$ is a locally Lipschitz vector field
\item $\mathbb{E}[ U_{n+1}^N \vert \mathcal{F}_n^N ]=0$
\item There exists $\Gamma \geqslant 0$ such that $\Vert U_n^N \Vert \leqslant \sqrt{\Gamma}$
\end{enumerate}

$F$ will be called the \textit{mean field} associated to $X^N$. The following results may be proved under a broader set of hypotheses but these ones will be sufficient for our purposes.

We denote by $\lbrace \varphi_t \rbrace$ the flow induced by $F$. So as to compare the trajectory of $\lbrace \varphi_t \rbrace$ with those of $(X_n^N)$ it's convenient to introduce the continuous in time process $\hat{X}^N : \R \to \R^m$ defined by
\[\hat{X}^N(k/N)=X_k^N \quad \forall k \in \N\]
and extended by linear interpolation on every $[k/N,(k+1)/N]$.

Let
\[D^N(T)=\max_{0\leqslant t \leqslant T} \Vert \hat{X}^N(t)-\varphi_t(X_0^N)\Vert\]
be the variable measuring the distance between the trajectories $t \mapsto \hat{X}^N(t)$ and $t \mapsto \varphi_t(X_0^N)$.

\begin{thm}
For every $T >0$, there exists $c >0$ (depending only on $F$,$\Gamma$ and $T$) such that, for every $\ep >0$, and for $N$ large enough :
\[\mathbb{P} [D^N(T) \geqslant \varepsilon ] \leqslant 2m\e{-\ep^2cN}\]
\end{thm}

See \cite{BenWei03} for a proof.

We now suppose that $X^N=(X^N_n)_{n\in \N}$ is a Markov chain defined on a countable set $K_N$. We also suppose that $X^N$ admits at least one invariant probability measure $\pi_N$.

\begin{thm}
The limit set of $\lbrace \pi_N \rbrace $ contains only probability measures that are invariant for the flow $\varphi_t$.
\end{thm}

See \cite{BenWei03} for a proof

%
%
%
%
%

For more information on the subject see e.g. \cite{BenWei03} and \cite{benaim-98}.
\addtocontents{toc}{\protect\setcounter{tocdepth}{-1}}
\section*{Acknowledgements}

The author wishes to thank Michel Bena\"im and Mathieu Faure for their advices and proofreading on this work.

\nocite{Kif88}
\nocite{FauSch11}

\bibliographystyle{plain}
\bibliography{these}

\end{document}